\documentclass[preprint, 12pt]{elsarticle}

\usepackage{amssymb}
\usepackage{amsmath}
\usepackage{algorithmic, multirow, color}
\def\bx{\mathbf{x}}
\def\bz{\mathbf{z}}

\def\bdelta{\boldsymbol{\delta}}

\journal{Computer Methods in Applied Mechanics and Engineering}

\begin{document}

\begin{frontmatter}



\title{A categorical embedding discontinuity-capturing shallow neural network for anisotropic elliptic interface problems}


\author[1]{Wei-Fan Hu}\ead{wfhu@math.ncu.edu.tw}

\author[2]{Te-Sheng Lin}\ead{teshenglin@nycu.edu.tw}
\author[3]{Yu-Hau Tseng\corref{cor1}}\ead{yhtseng@go.nuk.edu.tw}
\author[2]{Ming-Chih Lai}\ead{mclai@nycu.edu.tw}
\cortext[cor1]{Corresponding author}

\affiliation[1]{organization={Department of Mathematics, National Central University}, 
	postcode={32001}, city={Taoyuan}, country={Taiwan}}
\affiliation[2]{organization={Department of Applied Mathematics, National Yang Ming Chiao Tung University}, 
	postcode={30010}, city={Hsinchu}, country={Taiwan}}
\affiliation[3]{organization={Department of Applied Mathematics, National University of Kaohsiung}, 
	postcode={81148}, city={Kaohsiung}, country={Taiwan}}

%% Abstract
\begin{abstract}
In this paper, we propose a categorical embedding discontinuity-capturing shallow neural network for anisotropic elliptic interface problems. The architecture comprises three hidden layers: a discontinuity-capturing layer, which maps domain segments to disconnected sets in a higher-dimensional space; a categorical embedding layer, which reduces the high-dimensional information into low-dimensional features; and a fully connected layer, which models the continuous mapping. This design enables a single neural network to approximate piecewise smooth functions with high accuracy, even when the number of discontinuous pieces ranges from tens to hundreds. By automatically learning discontinuity embeddings, the proposed categorical embedding technique avoids the need for explicit domain labeling, providing a scalable, efficient, and mesh-free framework for approximating piecewise continuous solutions. To demonstrate its effectiveness, we apply the proposed method to solve anisotropic elliptic interface problems, training by minimizing the mean squared error loss of the governing system. Numerical experiments demonstrate that, despite its shallow and simple structure, the proposed method achieves accuracy and efficiency comparable to traditional grid-based numerical methods.
\end{abstract}


\end{frontmatter}



%%%%%%%%%%%%%%%%%%%%%%%%%%%%%%%%%%%%%%%%%
\section{Introduction}
%%%%%%%%%%%%%%%%%%%%%%%%%%%%%%%%%%%%%%%%%
Anisotropic elliptic interface problems are crucial to real-world applications, enabling accurate modeling and prediction in multiphase fluid flows, composite materials, biological systems, and geological processes. For example, the growth of liquid crystals involves interfaces between liquid and solid phases~\cite{WDB16, WDJH18}, so that the heat and mass transports during this process are anisotropic due to the underlying crystal structure. When modeling semiconductor heterostructures such as quantum wells, where abrupt changes in composite materials create interfaces, anisotropic elliptic partial differential equations can be used to simulate spin transport and diffusion~\cite{Tsymbal12} within these structures. Furthermore, placing anisotropic fluids such as nematic liquid crystals inside a Hele-Shaw cell~\cite{KSP98, FKSP01} enables the study of the interplay between the fluid-inherent anisotropy and the cell-confined geometry.

Owing to the essential role in many natural phenomena and industrial applications, finding effective numerical solutions to anisotropic elliptic interface problems is a significant task in the scientific computing community. High-order numerical methods have been developed for elliptic interface problems~\cite{J23, ZY24}; however, the case involving anisotropic diffusion tensors presents significantly greater challenges. The primary difficulty arises because the diffusion tensor is typically piecewise continuous, while the solution itself exhibits discontinuities across interfaces. To address this numerically, nonlinear finite volume schemes have been proposed to preserve solution positivity~\cite{LSSV07, ZLYS20}, alongside other approaches~\cite{PWXY22, XZ24, DFL22}. More recently, Li et al.~\cite{LY22} employed deep neural networks with a first-order formulation to approximate such problems. Their results only demonstrate smooth solutions; interface problems with piecewise continuous solutions remain considerably difficult. The low global regularity of these solutions makes the design of efficient and accurate numerical methods both essential and highly nontrivial.

Recent advances have demonstrated the effectiveness of neural network-based methods in tackling problems that remain challenging for classical numerical approaches, including inverse problems~\cite{PMAG21}, triangulated mesh prediction~\cite{GJF20}, and solutions with corner singularities~\cite{SF73}. From a theoretical perspective, the expressive power of neural networks has been rigorously established~\cite{Barron02, Cybenko89, HSW89, Pinkus99}, showing that continuous functions can be approximated to arbitrary precision. This provides a versatile function representation that can be leveraged to tackle a wide range of scientific and engineering problems. Notable examples include Physics-Informed Neural Networks (PINNs)~\cite{DP94, RPK19} and the Deep Ritz Method~\cite{EY18}, both of which have inspired extensive research on solving partial differential equations (PDEs) in complex geometries and high-dimensional settings. More recently, operator-learning frameworks such as the Deep Operator Networks~\cite{LJPZK21} and Fourier Neural Operator~\cite{LKALBSA20} have emerged, which directly approximate nonlinear operators mapping between function spaces, thereby expanding the scope and applicability of neural network approaches in scientific computing. On the other hand, for discontinuous functions, Llanas et al.~\cite{LLS08} proved that neural networks can approximate piecewise continuous functions almost uniformly. Their approach, however, relied on continuous activation functions to approximate step functions, which inevitably introduces an overall smoothness to the network representation. An alternative is to employ discontinuous activation functions~\cite{SL02, FN03, WM08}, commonly encountered in electronic circuits with switches or dry friction, which have shown effectiveness in optimization tasks such as linear and quadratic programming~\cite{LW08}. Nevertheless, how to systematically train such neural networks remains an open research question.

In this work, we propose to solve anisotropic elliptic interface problems within the PINNs framework by introducing a methodology for accurately representing discontinuous solutions composed of multiple smooth pieces. He et al.~\cite{HHM22} have attempted to partition the input space into subdomains separated by discontinuities, with individual neural networks assigned to each partition. While intuitive, this strategy suffers from scalability issues, as the number of required networks grows proportionally with the number of partitions. Alternative formulations based on weak form losses~\cite{LCLHL22} alleviate some of these difficulties but are hindered by the limited accuracy of Monte Carlo integration. To overcome these challenges, Hu et al. introduced the Discontinuity-Capturing Shallow Neural Network (DCSNN)~\cite{HLL22}, a single network architecture capable of accurately representing piecewise continuous functions. In this paper, we propose the categorical embedding DCSNN, a generalized extension of the original DCSNN, which enables the network to systematically maintain high accuracy even when approximating functions with a large number of discontinuous pieces, ranging from tens to hundreds. By automatically learning the embedding discontinuity information, the categorical embedding DCSNN avoids the need for explicit domain labeling and provides a scalable, efficient, and accurate framework for mimicking piecewise continuous solutions for anisotropic elliptic interface problems.

The rest of this paper is organized as follows. In Section~2, we introduce the mathematical model of anisotropic elliptic interface problems. In Section~3, we illustrate the network structures of DCSNN alongside the category-embedding structure for piecewise smooth functions, and apply the proposed neural network for anisotropic elliptic interface problems in Section~4. Section~5 provides a series of numerical experiments to demonstrate the capability and accuracy, and the concluding remarks are shown in Section~6.

%%%%%%%%%%%%%%%%%%%%%%%%%%%%%%%%%%%%%%%%%
\section{Anisotropic elliptic interface problems}\label{Sec:anisotropic}
\begin{figure}[ht]
	\centering
	\includegraphics[scale=0.35]{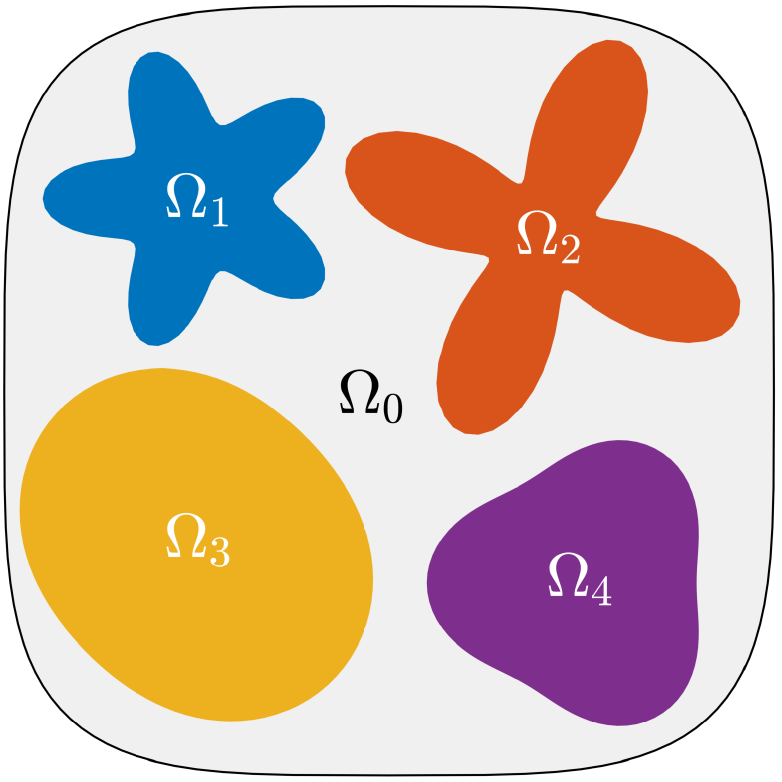}
	\caption{A two-dimensional domain, $\Omega = \cup_{\ell=0}^4\Omega_\ell$, contains four disjoint subdomains.}
	\label{Fig:domain_2d}
	\end{figure}

%%%%%%%%%%%%%%%%%%%%%%%%%%%%%%%%%%%%%%%%%
Let $\Omega\subseteq \mathbb{R}^d$ be a bounded, simply connected domain, together with $L$ disjoint subdomains, $\Omega_1,\Omega_2,\cdots,\Omega_L\subset \Omega$, and $\Omega_i \cap \Omega_j = \emptyset$ for $i\ne j$. Assume that the boundary of each subdomain, $\Omega_\ell$, denoted by $\Gamma_\ell$, is a $(d-1)$-dimensional closed $C^1$-hypersurface, and the complementary part of the union of all subdomains is represented by $\Omega_0$ so that $\Omega = \bigcup_{\ell=0}^L\Omega_\ell$. An illustrative two-dimensional example with $L=4$ is shown in Fig.~\ref{Fig:domain_2d}. In addition, in the rest of this paper, we assume that $u$ is a $d$-dimensional piecewise function $u:\Omega\rightarrow\mathbb{R}$ which is smooth within each subdomain $\Omega_\ell$ and exhibits jump discontinuities across each interface $\Gamma_\ell$.

The $d$-dimensional anisotropic elliptic interface problem with nonhomogeneous jump conditions is stated as follows:
\begin{subequations}\label{Eq:aniso}
\begin{align}
& \nabla\cdot(A(\bx)\nabla u(\bx)) - \lambda(\bx) u(\bx) = f(\bx) \quad \mbox{\;\;in\;\;}\Omega\setminus\bigcup\limits_{\ell=1}^L\Gamma_\ell,\\
& [u] = v_\ell(\bx), \quad [A\nabla u \cdot \mathbf{n}] = w_\ell(\bx) \quad \mbox{\;\;on\;\;} \Gamma_\ell, \mbox{ for\;\;} \ell = 1, 2, \cdots, L,
\end{align}
\end{subequations}
where $\nabla$ and $\nabla\cdot$ are respectively the gradient and divergence operators with respect to the spatial variable $\bx = (x_1,x_2,\cdots,x_d)$, $A(\bx)\in\mathbb{R}^{d\times d}$ is a positive-definite matrix, and $\lambda(\bx)$ is a nonnegative scalar function. The objective solution $u$ is typically piecewise smooth when the source term $f$, the matrix $A(\bx)$, and $\lambda(\bx)$ are piecewise-defined functions which are smooth within each subdomain $\Omega_{\ell}$, representing direction-dependent or piecewise properties due to the heterogeneous materials. The bracket $[\cdot]$ denotes the jump quantity for function values approaching from $\Omega_0$ minus the one from the subdomain $\Omega_\ell$; $\mathbf{n}$ is the unit outward normal vector (pointing from $\Omega_\ell$ toward $\Omega_0$) defined on the interface $\Gamma_\ell$. 

To close the system, a boundary condition must be specified along $\partial\Omega$. Throughout this paper, we assume that the solution is imposed by the Dirichlet boundary condition $u(\bx)|_{\partial\Omega} = g(\bx)$, while other types of boundary condition (Neumann or Robin type) can be implemented easily without changing the main ingredient of the proposed method (see the implementation in the following).

%%%%%%%%%%%%%%%%%%%%%%%%%%%%%%%%%%%%%%%%%
\section{Neural network architecture for piecewise smooth functions}
%%%%%%%%%%%%%%%%%%%%%%%%%%%%%%%%%%%%%%%%%
As discussed earlier, a standard neural network utilizing smooth activation functions yields inherently smooth representations, making it unable to capture jump discontinuities. To approximate piecewise-smooth functions, we revisit our previous work on DCSNN~\cite{HLL22}, which demonstrates the ability to model the discontinuous behavior using a single shallow neural network architecture. Building on this foundation, we introduce several classification models, particularly the categorical embedding technique, which automatically learns the embedding discontinuity information of piecewise-smooth functions with a large number of pieces.

%%%%%%%%%%%%%%%%%%%%%%%%%%%%%%%%%%%%%%%%%
%\subsection{Discontinuity capturing shallow neural networks with scalar encoding}
\subsection{Scalar encoding discontinuity-capturing shallow neural networks}
%%%%%%%%%%%%%%%%%%%%%%%%%%%%%%%%%%%%%%%%%
The core idea of DCSNN~\cite{HLL22} is to transform a $d$-dimensional piecewise-continuous function into a continuous function defined in a ($d+1$)-dimensional space. This is done by introducing a mapping that performs classification, assigning each subdomain to a point in a mutually distinct point set. For this purpose, we define a scalar categorical function $z: \mathbb{R}^d \to \mathbb{R}$ as follows:  
	\begin{equation}\label{Eq:DCSNN}
	z(\bx) = \gamma_\ell \quad \text{if } \bx \in \Omega_\ell,
	\end{equation}  
where $ \gamma_\ell $ for $\ell=0,1,\cdots,L$ are ``predefined'' constants. A natural choice is to set $\gamma_\ell = \ell$, whereby $\gamma_\ell$ serves as the region index and the label identifies the subdomain to which $\bx$ belongs. Alternatively, similar to target embedding categorization~\cite{B01}, one can incorporate input data information into the label, and define $\gamma_\ell$ based on function averages over each region: $\gamma_\ell = \bar{\gamma_\ell} = \int_{\Omega_\ell} u(\mathbf{x})\mathrm{d}\mathbf{x} / |\Omega_\ell|$. However, these values may be unavailable if $u$ itself is a function to be found (for instance, when $u$ is the solution of a PDE). %see Section~\ref{Sec:anisotropic} for illustration. 

We should point out that the selection of scalars, $\gamma_0,\gamma_1,\cdots,\gamma_L$, implicitly introduces an ordering for the subdomains, which significantly influences the training of subsequent neural networks, especially when dealing with a large number of disjoint subdomains. In addition, seeking the optimal constants a priori is generally impractical. Nevertheless, as long as those label values are distinct, the function $\mathbf{x}\to z(\mathbf{x})$ maps disjoint subdomains to a disconnected point set.

Using the above scalar categorical encoding function \( z(\mathbf{x}) \in \mathbb{R} \), we define an augmented function $U: \mathbb{R}^{d+1} \to \mathbb{R}$ that satisfies  
\begin{align}\label{Eq:U_DCSNN}
U(\bx, z(\bx)) = u(\bx), \quad \text{if } \bx \in \Omega.
\end{align}  
It is worth noting that $U$ is well defined and continuous on a disconnected subset of $\mathbb{R}^{d+1}$, namely $\{ (\bx, z(\bx)) \in \mathbb{R}^{d+1} \mid \bx \in \Omega \}$. Furthermore, this function $U$ can be continuously extended to the entire $\mathbb{R}^{d+1}$ as ensured by the Tietze extension theorem~\cite{Munkres}. The remaining task of DCSNN is to construct a shallow neural network to approximate the augmented continuous function $U$.

Leveraging the rich expressiveness of neural networks, the extension function $U$ can be approximated by a fully-connected one-hidden-layer neural network, $U_\mathcal{N}$, which takes the form
	\begin{align}
	U_\mathcal{N}(\bx, z(\bx);\Theta) = \sum_{j=1}^N c_j \sigma(W_j [\bx, z]^\top + b_j),
	\end{align}
where $N$ is the number of neurons in the hidden layer and $\sigma$ is the activation function. The notation $\Theta$ denotes a vector collecting all the trainable parameters, including the weights, $c_j \in \mathbb{R}$ and $W_j \in \mathbb{R}^{1 \times (d + 1)}$, and the biases, $b_j \in \mathbb{R}$. In addition, we define a neural function $u_\mathcal{N}(\bx)$ through $U_\mathcal{N}$, representing inherently piecewise continuous functions as
	\begin{align}
	u_\mathcal{N}(\bx) = U_\mathcal{N}(\bx, z(\bx);\Theta).
	\end{align}
With the same input, this succinct notation makes it easy to estimate the error between the target function $u(\bx)$ and the neural network function $u_\mathcal{N}(\bx)$.

\begin{figure}[h]
\centering
\includegraphics[width=0.63\textwidth]{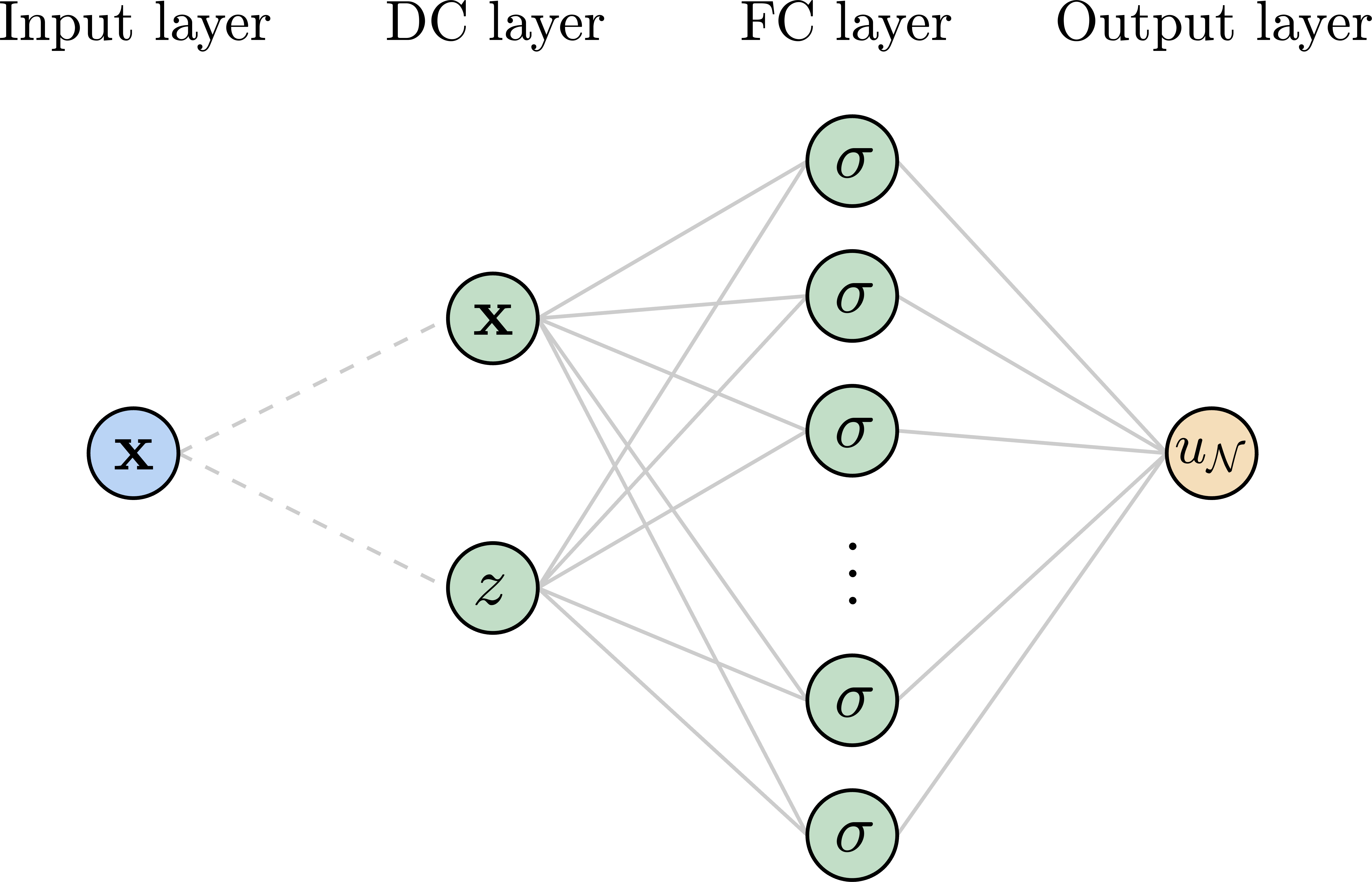}
\caption{Neural network structure of a scalar encoding DCSNN. The dashed/solid lines denote non-trainable/trainable parameters.}
\label{Fig:DCSNN_diagram}
\end{figure}

In summary, the scalar encoding DCSNN $u_\mathcal{N}$ consists of two layers, namely, a discontinuity-capturing (DC) layer and a fully-connected (FC) layer. The DC layer maps the input variable $\bx$ to $(\bx, z(\bx))$, which is pre-defined without training. On the contrary, the FC layer consists of $N$ neurons with trainable weights and biases. A schematic illustration of the network structure is shown in Fig.~\ref{Fig:DCSNN_diagram}. As a result, the scalar encoding DCSNN comprises a total of $N_p = (d+3)N$ trainable parameters.

%%%%%%%%%%%%%%%%%%%%%%%%%%%%%%%%%%%%%%%%%
\subsection{One-hot encoding neural networks}
%%%%%%%%%%%%%%%%%%%%%%%%%%%%%%%%%%%%%%%%%
Inspired by one-hot encoding in classification tasks, we propose an alternative representation of the categorical function. Specifically, we define a vector-valued function $\bz: \mathbb{R}^d \to \mathbb{R}^{L+1}$ to represent categorical variables in a binary format:
\begin{equation}\label{Eq:OH}
\bz(\bx) = \bdelta_{\ell+1}, \quad \text{if } \bx \in \Omega_{\ell},
\end{equation}
for $\ell = 0, 1, \cdots, L$, where $\bdelta_{\ell+1}$ is the $(l+1)$th vector of the canonical basis for $\mathbb{R}^{L+1}$. This mapping offers the advantage that the distance between the augmented features $\bz$ of any two subdomains remains constant, avoiding any implicit ordering in the augmentation coordinates. Mimicking the idea as introduced in the scalar encoding DCSNN, we then construct a $(d+L+1)$-dimensional continuous extension function $U$ that incorporates the vector-valued categorical function $\bz$ as follows:
\begin{align}\label{Eq:U_OH}
U(\bx, \bz(\bx)) = u(\bx), \quad \text{if } \bx \in \Omega.
\end{align}
This intermediate map $U$ can again be simply approximated by a shallow network $u_\mathcal{N}$ of the form
\begin{align}\label{Eq:UN_OH}
u_\mathcal{N}(\bx) = U_\mathcal{N}(\bx, \bz(\bx);\Theta) = \sum_{j=1}^N c_j \sigma(W_j [\bx, \bz]^\top + b_j),
\end{align}
where $\Theta$ represents a vector collecting all the trainable parameters. Notably, the usage of the vector-valued functions as augmented inputs has also been explored in~\cite{HCCCY23}, where it was employed to classify two subdomains within their network architecture.

Compared to the scalar encoding model~(\ref{Eq:DCSNN}), adopting the one-hot encoding model seems ideal in the absence of ordinal or nominal labels for each subdomain. However, approximating this continuous function $U$ using a neural network expression with the augmented input $(\bx, \bz) \in \mathbb{R}^{d+L+1}$ poses computational challenges as the number of subdomains increases (i.e., larger $L$), leading to a significant increase in the number of trainable parameters. Consequently, the number of unknown parameters of the network expression~(\ref{Eq:UN_OH}) now becomes $N_p = (d+L+3)N$, making the training process computationally intensive.

%%%%%%%%%%%%%%%%%%%%%%%%%%%%%%%%%%%%%%%%%
\subsection{Categorical embedding discontinuity-capturing shallow neural networks}
%%%%%%%%%%%%%%%%%%%%%%%%%%%%%%%%%%%%%%%%%
From the two models discussed above, it is evident that multiple options exist for defining the categorical function. The role of this function is to construct a mapping that both lifts each smooth component of the piecewise function into a higher-dimensional space and ensures their separation in that space. Consequently, there are infinitely many possible choices for such a mapping.

In this context, we propose a categorical embedding DCSNN, in which the objective is to \emph{learn} an optimal categorical function. This function aims to effectively capture the intrinsic properties (or features) within the function profiles of each subdomain (or category), and thus map similar categories closer together in a specific low-dimensional space. To this end, we define the categorical function $\bz: \mathbb{R}^{d} \to \mathbb{R}^{D}$ via a linear map:
\begin{align}\label{Eq:CE}
\bz(\bx) = E\bdelta(\bx),
\end{align}
where $\bdelta(\bx) = \bdelta_{\ell+1} \in \mathbb{R}^{L+1}$ for $\bx \in \Omega_\ell$, and $E \in \mathbb{R}^{D \times (L+1)}$ is the embedding matrix that classifies the $(L+1)$ subdomains into a low-dimensional embedded space of the selected dimension $D$ (with $D \leq L+1$). We remark that the categorical embedding function~(\ref{Eq:CE}) used here follows the same principle as the entity embedding model~\cite{GB16}, where high-cardinality categorical variables are embedded into low-dimensional Euclidean spaces.

Again, we look for an extension function $U:\mathbb{R}^{d+D} \to \mathbb{R}$ with the proposed categorical embedding map~(\ref{Eq:CE}) that satisfies
\begin{align}\label{Eq:U_CE}
U(\bx, \bz(\bx)) = U(\bx, E\bdelta(\bx)) = u(\bx), \quad \text{if } \bx \in \Omega.
\end{align}
At this stage, the embedding matrix $E$ is treated as an unknown weight matrix and can be learned through standard machine learning optimization techniques. Additionally, we emphasize that the importance of this categorical step lies in mapping low-dimensional discontinuous functions into high-dimensional smooth functions.

As before, it involves approximating the extension function $U$ by constructing a categorical embedding DCSNN $u_\mathcal{N}$ through the intermediate map $U_\mathcal{N}$ that represents inherently piecewise functions as follows:
\begin{align}\label{Eq:UNN_CE}
u_\mathcal{N}(\bx) = U_\mathcal{N}(\bx, E\bdelta(\bx);\Theta) = \sum_{j=1}^N c_j \sigma(W_j [\bx, E\bdelta(\bx)]^\top + b_j),
\end{align}
where $\Theta$ represents a vector that collects all the trainable parameters, including those in the embedding feature matrix $E$. The number of total trainable parameters, including weights, biases, and embedding matrix $E$, is $N_p = (d+D+2)N+(L+1)D$.

\begin{figure}[h]
\centering
\includegraphics[width=0.7\textwidth]{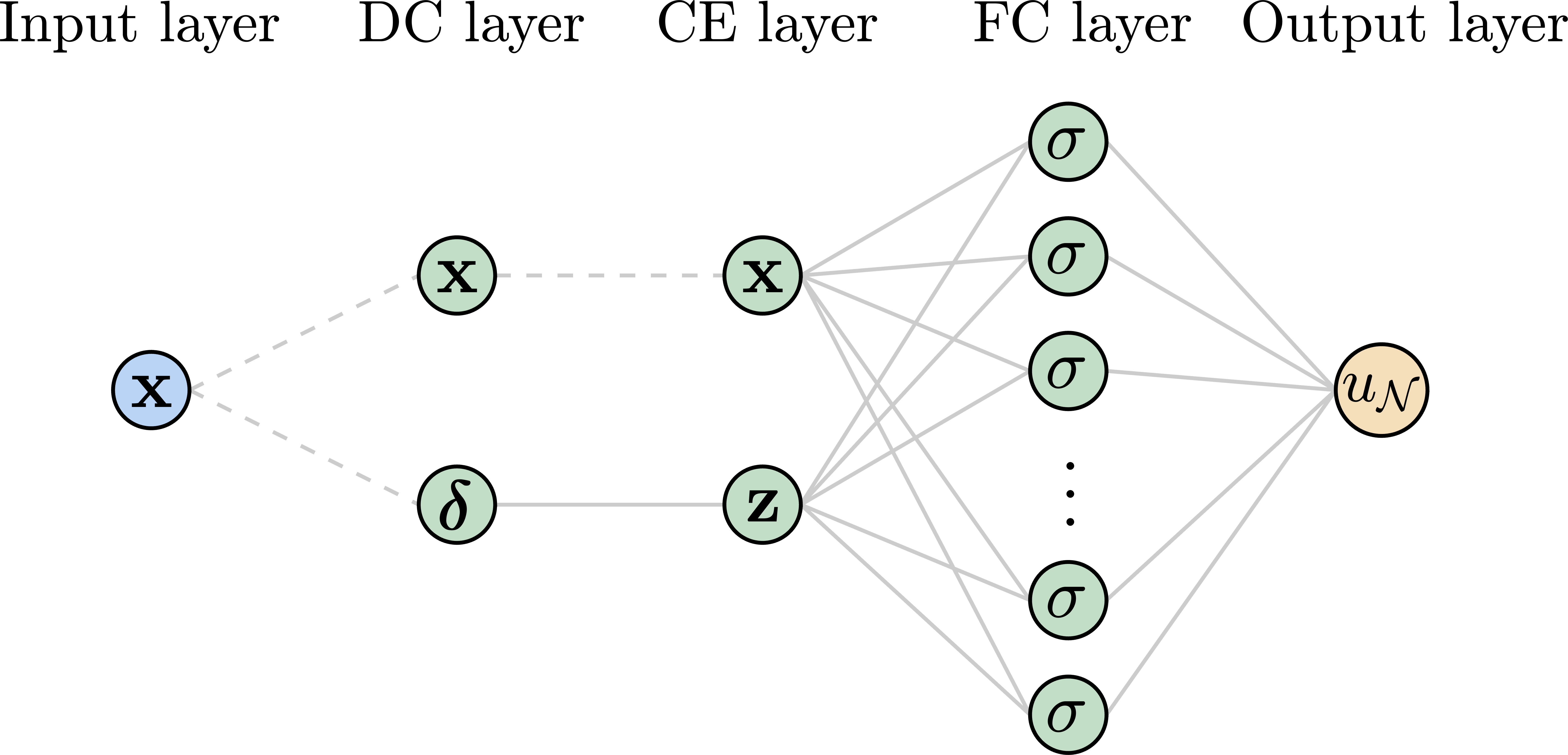}
\caption{Neural network structure of categorical embedding DCSNN. The dashed/solid lines denote non-trainable/trainable parameters.}
\label{Fig:NN_diagram}
\end{figure}

The categorical embedding DCSNN consists of three layers: a DC layer, a categorical embedding (CE) layer, and an FC layer. The DC layer maps the input variable $\bx$ to $(\bx, \bdelta(\bx))$, which is predefined and does not require training. In contrast, the CE layer learns the embedding (or weight) matrix $E$ via the mapping $\bz(\bx) = E\bdelta(\bx)$. The FC layer consists of $N$ neurons with trainable weights and biases. A schematic representation of the network structure is shown in Fig.~\ref{Fig:NN_diagram}. In this paper, we focus on the simple structure defined in Eq.~(\ref{Eq:UNN_CE}). While replacing the FC layer with a deep neural network or adding residual connections to the proposed model is straightforward.  
%Our model can easily be extended by deepening the FC layer with more hidden layers or by adding residual connections, but for this work, we focus on the simpler structure defined in~(\ref{Eq:UNN_CE}).

We provide several remarks on the proposed model in the following. Firstly, the categorical embedding technique generalizes both scalar encoding and one-hot encoding models depending on the choice of the embedding matrix $E$. Particularly, setting $E = [\gamma_0, \gamma_1, \cdots, \gamma_{L}]$ recovers the scalar encoding model (see Eq.~(\ref{Eq:DCSNN})), while choosing the identity matrix $E = I_{L+1}$ reverts to the one-hot encoding approach (see Eq.~(\ref{Eq:OH})). Secondly, one can also generalize the proposed linear embedding by including nonlinearity. For instance, applying an activation function to the categorical embedding, i.e., setting $\bz(\bx) = \sigma(E\bdelta(\bx))$, may enhance the classification of the encoding labels in the embedded space. Thirdly, selecting the best low-dimensional embedding, or the optimal reduced dimension $D$, of the neural network function $u_\mathcal{N}$ in Eq.~(\ref{Eq:UNN_CE}) remains a subject of ongoing research. To provide insight into the choice of $D$, we will present a systematic study of function approximation experiments in Section~\ref{subsec:test_func_approx}. Lastly, the existence of a smooth extension function $U$ on the entire domain is again guaranteed by the Tietze extension theorem~\cite{Munkres}, ensuring that such a function can be approximated using neural networks for the above three models. 

%%%%%%%%%%%%%%%%%%%%%%%%%%%%%%%%%%%%%%%%%
\subsection{Training method}
%%%%%%%%%%%%%%%%%%%%%%%%%%%%%%%%%%%%%%%%%
In the training procedure of the network function $u_\mathcal{N}(\bx)$ in Eq.~(\ref{Eq:UNN_CE}), we utilize the commonly used mean-squared error loss function. As in supervised learning tasks, we approximate a piecewise-defined function $u(\bx)$ using the neural network $u_\mathcal{N}(\bx)$. We choose a set of training data, $\{(\bx^i, u^i)\}_{i=1}^M$, where $\bx^i\in \Omega$ and $u^i=u(\bx^i)$, and specify the embedding function $\bz(\bx)$. The loss function under the supervised learning framework is then defined as
\begin{align}\label{Eq:loss_approx}
\mbox{Loss}(\Theta) = \frac{1}{M}\sum_{i=1}^M\left( u^i - u_\mathcal{N}(\bx^i) \right)^2.
\end{align}
Recall that $u_\mathcal{N}\left(\bx^i\right)=U_\mathcal{N}\left(\bx^i,\bz\left(\bx^i\right);\Theta\right)$, and $\Theta$ is a vector collecting all the trainable parameters. 

Throughout the numerical experiments in this paper, we employ the Levenberg-Marquardt (LM) method~\cite{More78}, which is particularly effective due to the adopted least-squares loss formulation. The update at step $(k+1)$ in the LM method is expressed as  
\begin{align}\label{Eq:LM}
\Theta^{(k+1)} = \Theta^{(k)} + (J^\top J + \mu I)^{-1}\left[ J^\top\left(\mathbf{u}-\mathbf{u}_\mathcal{N}\left(\Theta^{(k)}\right)\right) \right],  
\end{align}  
where $\mu > 0$ is a tunable damping parameter; \(\mathbf{u}\) and \(\mathbf{u}_\mathcal{N}\left(\Theta^{(k)}\right)\) denote the vectors collecting the data \(u^i\) and \(U_\mathcal{N}\left(\bx^i,\bz\left(\bx^i\right);\Theta^{(k)}\right)\) in the loss function~(\ref{Eq:loss_approx}), respectively. Here, \(J\) is the Jacobian matrix defined as \(\partial \mathbf{u}_\mathcal{N}\left(\Theta^{(k)}\right) / \partial \Theta\). Obviously, the primary computational expense in the LM update step \( \Delta \Theta \) (i.e., the matrix-vector multiplication in the right-hand side of Eq.~(\ref{Eq:LM})) arises from solving the regularized least-squares problem with the damping parameter \( \mu \). This step is typically performed using Cholesky factorization. However, when \( \mu \ll 1 \), the condition number of \( J^\top J + \mu I \) can become extremely large, leading to instability and poor approximation of the update step. To mitigate this issue, we can compute the parameter update \( \Delta \Theta \) using QR factorization by solving the following linear system:  
\begin{align*}
\left[
\begin{array}{c}
J\\
\sqrt{\mu}I
\end{array}\right]\Delta\Theta = 
\left[
\begin{array}{c}
\mathbf{u}-\mathbf{u}_\mathcal{N}\left(\Theta^{(k)}\right)\\
\mathbf{0}
\end{array}\right].
\end{align*}
Note that the network model with the loss function~(\ref{Eq:loss_approx}) is conventionally trained using optimizers such as Adam~\cite{HS18} or L-BFGS~\cite{LN89}, which are widely adopted in the literature. In subsequent experiments, we will compare the efficiency of various approaches with these optimizers during the training process.

%%%%%%%%%%%%%%%%%%%%%%%%%%%%%%%%%%%%%%%%%
\section{Categorical embedding shallow neural network for anisotropic elliptic interface problems}\label{sec:pinn_aeips}
%%%%%%%%%%%%%%%%%%%%%%%%%%%%%%%%%%%%%%%%%
%In this section, we apply the proposed categorical embedding model~(\ref{Eq:UNN_CE}) as a solution representation to address the anisotropic elliptic interface problem~(\ref{Eq:aniso}). The derivation of jump discontinuity and derivatives of the categorical embedding function $u_\mathcal{N}(\bx)$ in Eq.~(\ref{Eq:UNN_CE}) are presented in Section~\ref{subsec:property}. Section~\ref{subsec:pinn_aeips} employs a PINN framework for solving anisotropic elliptic interface problems.

In this section, we apply the proposed categorical embedding model~(\ref{Eq:UNN_CE}) as a solution representation for the anisotropic elliptic interface problem~(\ref{Eq:aniso}). The derivation of the jump discontinuities and the derivatives of the categorical embedding function $u_\mathcal{N}(\bx)$ in Eq.~(\ref{Eq:UNN_CE}) is presented in Section~\ref{subsec:property}, while Section~\ref{subsec:pinn_aeips} introduces a PINN framework for solving such a problem.
%%%%%%%%%%%%%%%%%%%%%%%%%%%%%%%%%%%%%%%%%
\subsection{Discontinuities and derivatives}\label{subsec:property}
% better move this subsection to neural networks for AEIPs
%%%%%%%%%%%%%%%%%%%%%%%%%%%%%%%%%%%%%%%%%

\paragraph{Values at jump discontinuities}

Generally, the function values of a discontinuous function at a breakpoint hold less significance. Of greater importance is the jump value at the breakpoint, which represents the difference between the limiting values approaching from two opposite directions. For instance, denoting the boundary of subdomain $\Omega_1$ as $\Gamma_1$, the jump quantity for the function value at $\bx\in\Gamma_1$ is given by
\begin{equation}
[u]_{\bx} = \lim_{\bx^+\to\bx}u(\bx^+) - \lim_{\bx^-\to\bx}u(\bx^-),
\end{equation}
where $\bx^+\in\Omega_0$ and $\bx^-\in\Omega_{1}$ (see Fig.~\ref{Fig:domain_2d} for example). Surprisingly, this information naturally emerges within the intermediate map $U$ from Eq.~(\ref{Eq:U_CE}). That is, taking the categorical function $\bz$ in Eq.~(\ref{Eq:CE}), the limiting value approaching from $\Omega_0$ follows
\begin{equation}
\lim_{\bx^+\to\bx}u(\bx^+) = \lim_{\bx^+\to\bx} U(\bx^+, \bz(\bx^+)) = U(\bx, E\bdelta_0),
\end{equation}
while the other limit is $\lim_{\bx^-\to\bx}u(\bx^-) = U(\bx, E\bdelta_1)$. Thus the jump value at $\bx\in\Gamma_1$, which typically requires taking one-sided limits of the function $u$, can be calculated directly by evaluating $[u]_{\bx} = U(\bx,E\bdelta_1) - U(\bx,E\bdelta_0)$. The same manner applies to the jump of the network function $u_\mathcal{N}$ as
\begin{equation}\label{Eq:1d_jump}
[u_\mathcal{N}]_{\bx} = U_\mathcal{N}(\bx, E\bdelta_0;\Theta) - U_\mathcal{N}(\bx, E\bdelta_1;\Theta).
\end{equation}
The rationale behind this outcome is the smooth nature of the intermediate map $U$. Consequently, determining the jump quantity can be accomplished easily through the function evaluation of the intermediate map, eliminating the need to take limits of the function.

\paragraph{Derivative evaluation}
%There is also a simple relation between the derivative of the network function $u_\mathcal{N}$ and the intermediate map $U_\mathcal{N}$. 
Except at breakpoints, the categorical function $\bz(\bx)$ remains constant everywhere; hence, its derivative is zero, implying that the Jacobian matrix $\nabla\bz = \mathbf{0}$. To express the gradient of network function $u_\mathcal{N}$ in terms of the intermediate map $U_\mathcal{N}$, the chain rule gives
\begin{equation}\label{Eq:1d_der}
\nabla u_\mathcal{N} = \nabla_{\bx} U_\mathcal{N} + \nabla\bz\nabla_{\bz}U_\mathcal{N}  = \nabla_{\bx} U_\mathcal{N},
\end{equation}
where $\nabla$ is the conventional gradient operator; $\nabla_{\bx}$ and $\nabla_{\bz}$ denote differentiating only with respect to $\bx$ and $\bz$, respectively. As a result, the partial derivatives of $u_\mathcal{N}$ can be calculated equivalently by taking the partial derivatives of $U_\mathcal{N}$ with respect to the original variables $\bx$. Higher-order derivatives of $u_\mathcal{N}$ are calculated in the same manner.

\subsection{PINN framework}\label{subsec:pinn_aeips}
We now apply the methodology of physics-informed learning machinery~\cite{RPK19} for solving Eq.~(\ref{Eq:aniso}) as follows. To find the network parameters in the expression~(\ref{Eq:UNN_CE}), we convert Eq.~(\ref{Eq:aniso}) to an optimization problem via a loss function. Namely, we first choose the sets of training points in the domain, on the domain boundary, and along all the interfaces, as
\begin{equation*}
\left\{\bx^i\right\}_{i=1}^M\subseteq\Omega, \qquad
\left\{\bx_{\partial\Omega}^j\right\}_{j=1}^{M_b}\subseteq\partial \Omega, \qquad
\left\{\bx_{\Gamma_\ell}^k\right\}_{k=1}^{M_{\Gamma_\ell}}\subseteq\Gamma_\ell,
\end{equation*}
respectively. The loss function is then defined as 
\begin{eqnarray}\label{Eq:loss_aniso}
\mbox{Loss}(\Theta) &=& \frac{1}{M}\sum_{i=1}^M\left( \nabla\cdot(A(\bx^i)\nabla u_\mathcal{N}(\bx^i)) - \lambda(\bx^i) u_\mathcal{N}(\bx^i) - f(\bx^i) \right)^2 \nonumber\\
& +& \frac{1}{M_b}\sum_{j=1}^{M_b} \left( u_\mathcal{N}(\bx_{\partial\Omega}^j) - g(\bx_{\partial\Omega}^j) \right)^2 \\
& +& \sum_{\ell=1}^L\frac{1}{M_{\Gamma_\ell}}\left( \sum_{k=1}^{M_{\Gamma_\ell}} \left( [u_\mathcal{N}]-v_\ell(\bx_{\Gamma_\ell}^k) \right)^2 +  \left( [A\nabla u_\mathcal{N} \cdot \mathbf{n}]-w_\ell(\bx_{\Gamma_\ell}^k) \right)^2 \right), \nonumber
\end{eqnarray}
which consists of the mean squared residual for each equation in the original problem~(\ref{Eq:aniso}), following the same principle as the PINN-type loss~\cite{DP94, RPK19}. We also recall that $\Theta$ represents the set of trainable parameters, and the objective is to find $\Theta$ that minimizes the loss function. As each term in the loss again takes the form of least-squares errors, we can train the model efficiently using the LM algorithm. 

As mentioned in Sec.~\ref{subsec:property}, we recall that the categorical function $\bz(\bx)$ is a piecewise constant vector function that has zero derivative over the interior of each subdomain, so that computing derivatives of $u_\mathcal{N}$, such as gradient or divergence, can be straightforwardly done using chain rule without any difficulty. For example, recall the relation~(\ref{Eq:1d_der}), we have
\begin{equation}
\nabla u_\mathcal{N}(\bx) = \nabla_{\bx} U_\mathcal{N}(\bx, \bz;\Theta), \quad \bx\in\Omega\setminus\bigcup\limits_{\ell=1}^L\Gamma_\ell,
\end{equation}
where $U_\mathcal{N}$ is the intermediate map and $\nabla_{\bx}$ is the gradient operator with respect to the $\bx$ variable only. Therefore, all the derivatives in Eq.~(\ref{Eq:loss_aniso}) are well-defined. As regards the jump conditions in Eq.~(\ref{Eq:loss_aniso}), they can be computed easily through function evaluations of $U_\mathcal{N}$ (see Sec.~\ref{subsec:property}). Similarly, the computation for the flux jump condition $[A\nabla u_\mathcal{N} \cdot \mathbf{n}]$ can be evaluated in the same manner as $[u_\mathcal{N}]$.

Moreover, the derivative terms involved in the loss model are commonly computed using automatic differentiation. However, in fact, thanks to the design of the shallow network structure~(\ref{Eq:UNN_CE}), one can easily derive the explicit form of the derivatives, which is much more efficient in practice. For example, the partial derivative with respect to the $k$-th spatial component is obtained by
\begin{align}
\frac{\partial u_\mathcal{N}}{\partial x_k}(\bx) = \frac{\partial U_\mathcal{N}}{\partial x_k}(\bx,\bz;\Theta) = \sum_{j=1}^N c_jW_{jk}\sigma'(W_j[\bx,\bz]^\top+b_j),
\end{align}
where $W_{jk}$ is the $k$-th component of the vector $W_j$; the prime notation of $\sigma$ means the derivative of the activation function. We mention that both the higher order or mixed partial derivatives to the target function and the Jacobian matrix (collecting all partial derivatives with respect to the trainable parameters $\Theta$ for each residual loss in Eq.~(\ref{Eq:loss_aniso})) involved in the LM training iteration can be implemented straightforwardly using the simple formulation (\ref{Eq:UNN_CE}). We also point out that since a simple structure with a moderate number of neurons is employed in the present network, the computational complexity and learning workload can be significantly reduced.

%%%%%%%%%%%%%%%%%%%%%%%%%%%%%%%%%%%%%%%%%
\section{Numerical Results}
%%%%%%%%%%%%%%%%%%%%%%%%%%%%%%%%%%%%%%%%%
This section presents a series of numerical experiments to demonstrate the capability of the proposed categorical embedding techniques for representing piecewise smooth functions and solving anisotropic elliptic interface problems. %\textcolor{red}{To reproduce the numerical experiments in this study, the code will be made publicly available at \texttt{https://github.com/teshenglin/DCNN\underline{ }CE}  after the paper is accepted. (NECESSARY?)}

%%%%%%%%%%%%%%%%%%%%%%%%%%%%%%%%%%%%%%%%%
\subsection{Examples for the encoding methods}\label{subsec:test_func_approx}
%%%%%%%%%%%%%%%%%%%%%%%%%%%%%%%%%%%%%%%%%
We present two examples to illustrate the capability of the proposed CE model, originally introduced as the categorical embedding technique, in approximating piecewise smooth functions. Additionally, we compare its performance with scalar encoding (SE) and one-hot encoding (OH) models. In all network models, we use $N = 50$ neurons in the FC layer, maintaining the same number of basis network functions across different approaches. Each neuron is employed with the sigmoid activation function. For each example, we randomly sample test points approximately 10 times the number of training points to evaluate the average $L^2$ and $L^\infty$ errors over 10 trial runs. The LM optimizer stops either after 1000 training steps or reaching the tolerance $\varepsilon = 10^{-15}$.

%%%%%%%%%%%%%%%%%%%%%%%%%%%%%%%%%%%%%%%%%
\subsubsection*{\textbf{Example 1}}
%%%%%%%%%%%%%%%%%%%%%%%%%%%%%%%%%%%%%%%%%
We consider a domain $\Omega\subset\mathbb{R}^2$ that is enclosed by the superellipse, $x_1^4+x_2^4 = 1$, which encapsulates four subdomains whose boundaries (or interfaces) are described by the polar curves as $r_1(\theta) = 0.3 - 0.1\cos(5\theta), r_2(\theta) = 0.35 - 0.2\sin(4\theta), r_3(\theta) = 0.45 - 0.05\sin(2\theta), r_4(\theta) = 0.35 - 0.05\cos(3\theta)$ with the center located at $(-0.5,0.5)$, $(0.4,0.4)$, $(-0.5,-0.4)$, $(0.5,-0.5)$, respectively. The domain is depicted in Fig.~\ref{Fig:domain_2d}. The target function $u$ is chosen as
\begin{equation}\label{Eq:function_approx}
u(x_1,x_2) =
\left\{
\begin{array}{ll}
\sin (x_1)\sin (x_2) & \mbox{\;\;if\;\;} (x_1,x_2) \in \Omega_0,\\
\exp(x_1-x_2)      & \mbox{\;\;if\;\;} (x_1,x_2) \in \Omega_1, \\
\cos(x_1+x_2)      & \mbox{\;\;if\;\;} (x_1,x_2) \in \Omega_2, \\
0.5\cosh(x_1+x_2)      & \mbox{\;\;if\;\;} (x_1,x_2) \in \Omega_3, \\
\ln(x_1+x_2+3)      & \mbox{\;\;if\;\;} (x_1,x_2) \in \Omega_4.
\end{array}\right.
\end{equation}
We minimize the loss function~(\ref{Eq:loss_approx}) using $M=1000$ randomly sampled training points, consisting of 880 points inside the domain $\Omega$ and 120 points along the domain boundary $\partial\Omega$. 

To demonstrate the training efficiency, we compare the performance of various optimizers discussed previously with the three categorical network models. For the scalar encoding model, we particularly set the nominal label values as $\gamma_\ell = \ell$ for $\ell = 0, 1, \cdots, 4$. 

\begin{figure*}[h]
\centering
\includegraphics[scale=0.23]{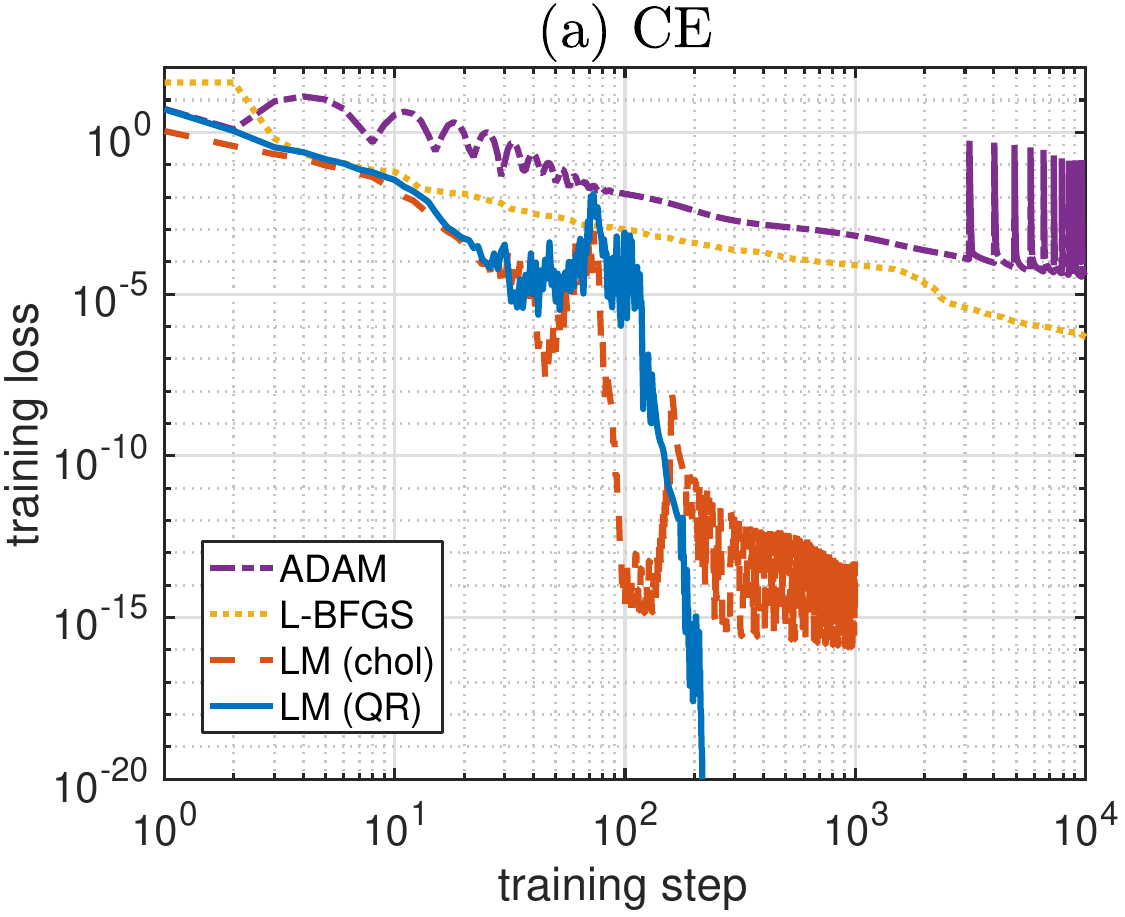}
\includegraphics[scale=0.23]{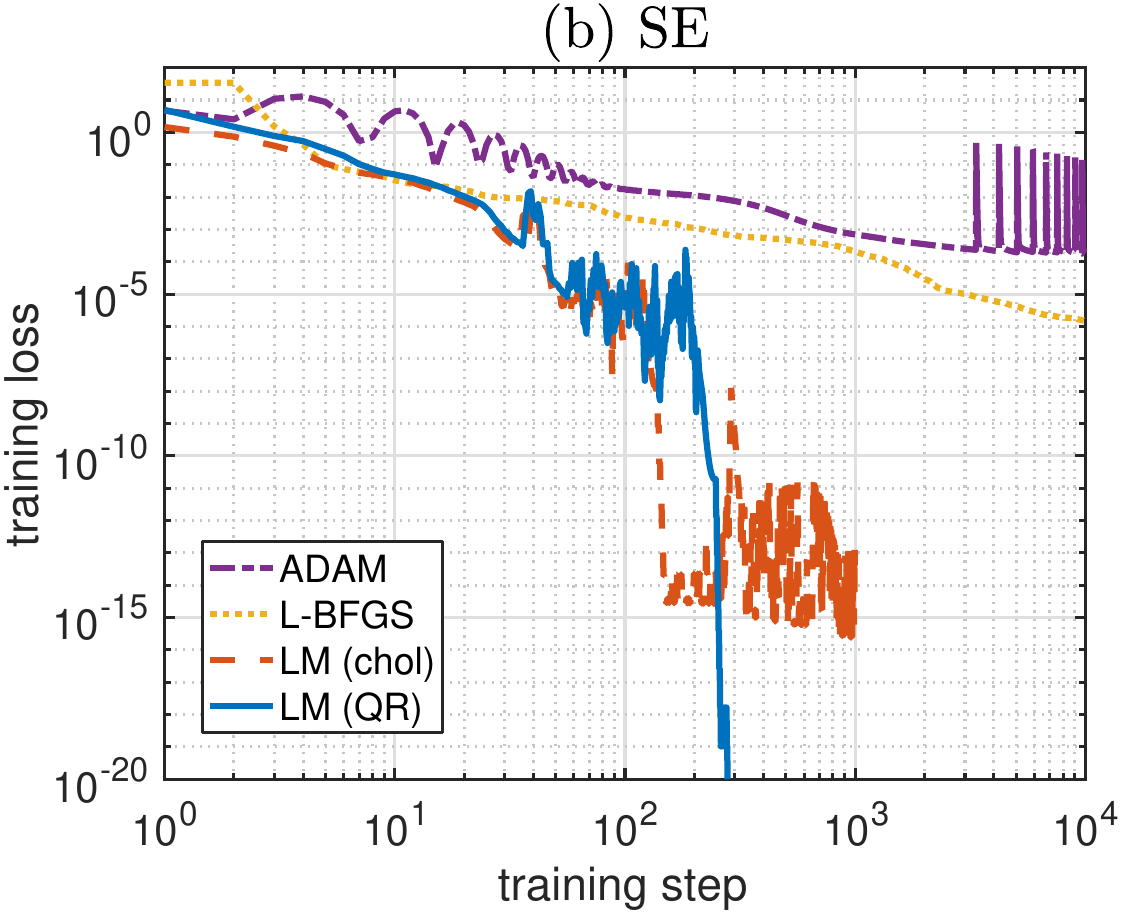}
\includegraphics[scale=0.23]{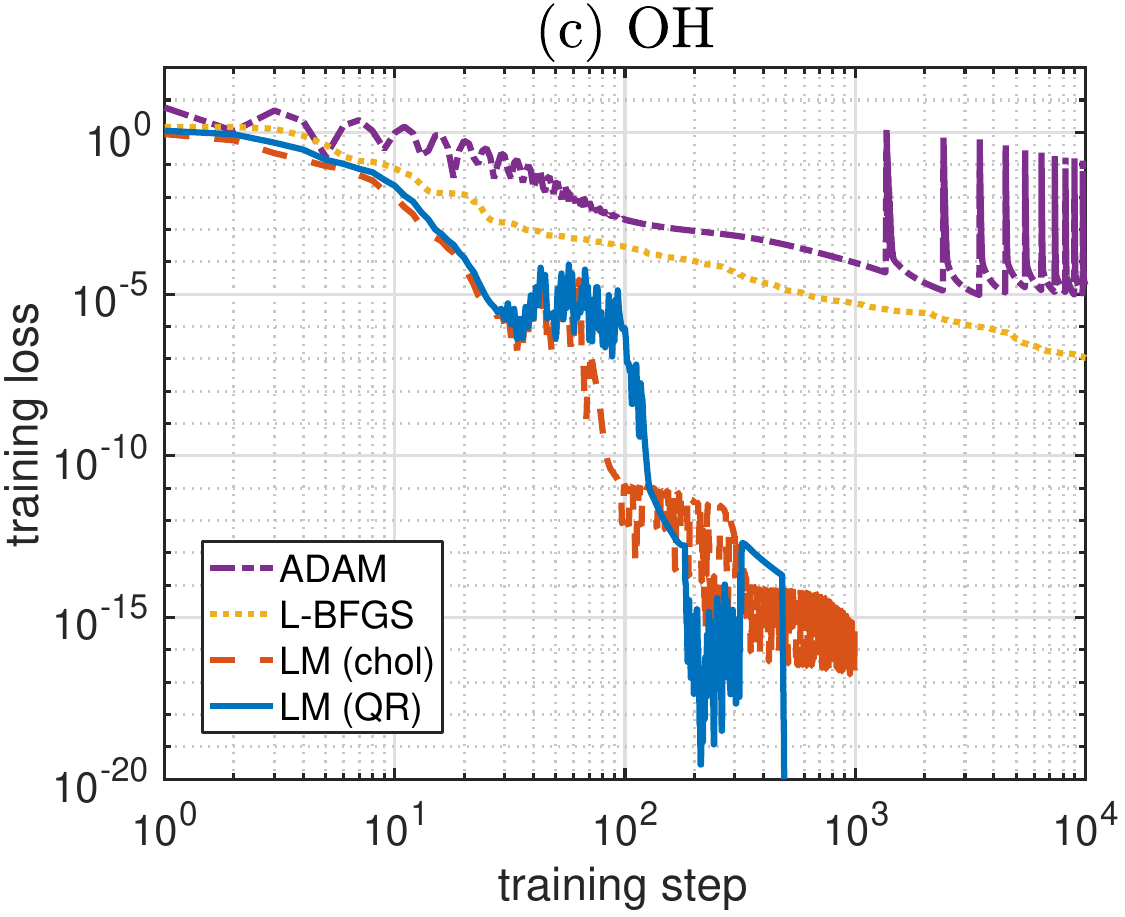}
\caption{Training history for (a) categorical embedding, (b) scalar encoding, and (c) one-hot encoding model with different optimizers.}
\label{Fig:optimizer}
\end{figure*}

The results are reported in Fig.~\ref{Fig:optimizer}, where all models exhibit a similar trend. Specifically, the gradient-based Adam optimizer (purple dash-dotted line) achieves a loss magnitude of approximately $10^{-5}$ but gets stuck in a local minimum, even after $10000$ training steps. Meanwhile, the quasi-Newton L-BFGS algorithm (yellow dotted line) performs slightly better than Adam but also nearly stagnates, resulting in a slow decrease of the loss value in subsequent steps. In contrast, the LM update strategy using Cholesky decomposition (red dashed line) reaches a local minimum as low as $10^{-15}$. As expected, the QR factorization method (blue solid line) achieves even lower loss values, around $10^{-20}$, within just a few hundred training steps.

Next, we demonstrate the capability of each network model for approximating the 2D discontinuous function.
Using the LM optimizer, the comparison results for all models are summarized in Table~\ref{Table:2d_multi_piece_comparison}. For the categorical embedding model, we select the reduced dimension $D = 1$. As anticipated, the one-hot encoding model, which requires the most parameters to be trained, achieves the most accurate prediction. However, the categorical embedding model, which has approximately half the number of parameters as the one-hot encoding model, yields nearly identical results. Furthermore, although both categorical embedding (with $D = 1$) and scalar encoding models use scalars to categorize each function piece, the categorical embedding model is certainly expected to perform better as it learns the optimal classification map in the embedding space.

\begin{table}[h]
\caption{Numerical results for approximating multi-piece function~(\ref{Eq:function_approx}). CE: categorical embedding model; SE: scalar encoding model; OH: one-hot encoding model.}
\label{Table:2d_multi_piece_comparison}
\centering
\begin{tabular}{c|ccc}
\hline
Method & $N_p$ & $L^2$ error & $L^\infty$ error\\
\hline
CE ($D = 1$) & 255 & 6.47E$-$08 & 2.09E$-$06\\
SE ($\gamma_\ell = \ell$)         & 250 & 1.10E$-$07 & 4.07E$-$06\\
OH                & 450 & 4.22E$-$08 & 1.06E$-$06\\
\hline
\end{tabular}
\end{table}

\begin{figure*}[!]
\centering
\includegraphics[scale=0.45]{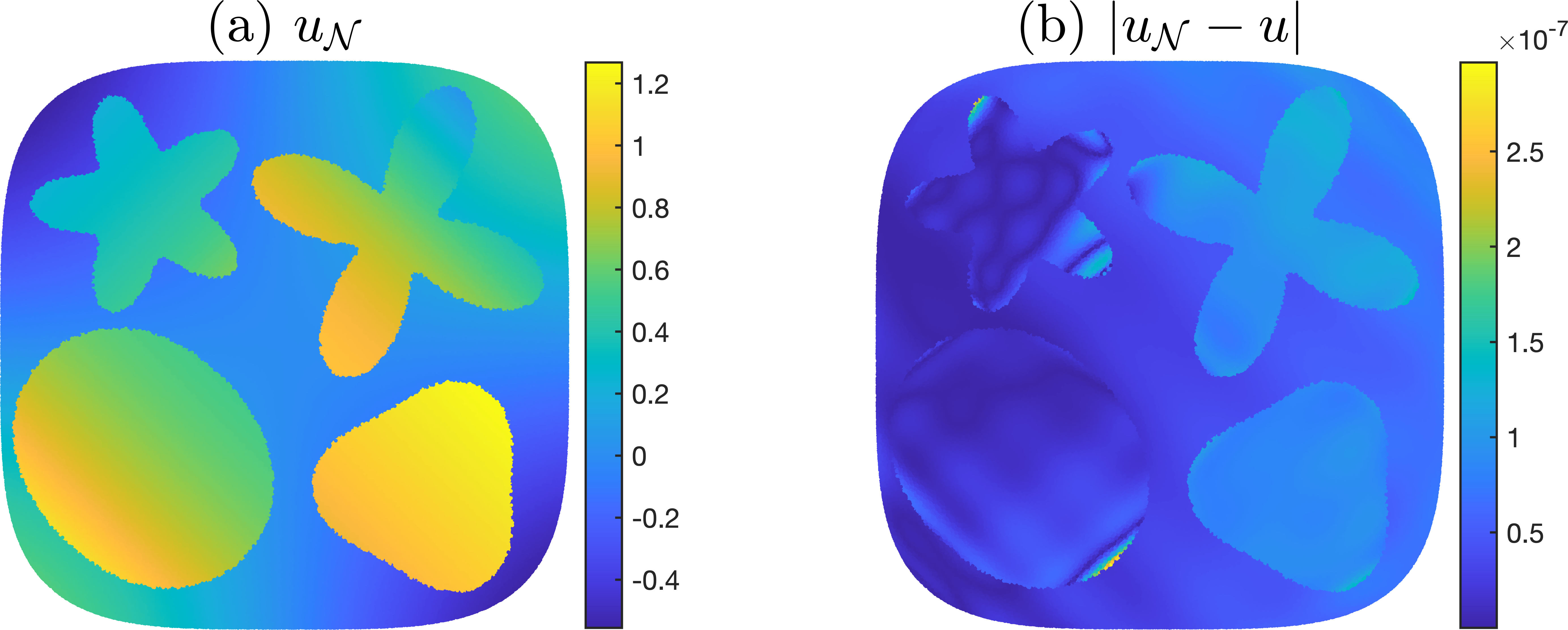}
\caption{(a) The trained CE model $u_\mathcal{N}$ for the piecewise-defined function in Example~1. (b) The pointwise absolute error $|u_\mathcal{N}-u|$. The maximum error is $\|u_\mathcal{N}-u\|_\infty = 2.96\times10^{-7}$.}
\label{Fig:func_mp}
\end{figure*}

Additionally, we present the network profile of the learned categorical embedding function in Fig.~\ref{Fig:func_mp}(a). As shown, the network model accurately captures all jump discontinuities sharply and represents the function well, with the pointwise absolute error, depicted in Fig.~\ref{Fig:func_mp}(b), being on the order of $10^{-7}$. It is observed that the significant errors mainly occur near the interfaces of each subdomain. This is because we only randomly sample the training points within the domain, but without additional information along the interfaces. 

%%%%%%%%%%%%%%%%%%%%%%%%%%%%%%%%%%%%%%%%%
\subsubsection*{\textbf{Example 2}}
%%%%%%%%%%%%%%%%%%%%%%%%%%%%%%%%%%%%%%%%%

In this example, we demonstrate the expressive power of the proposed CE network by approximating a discontinuous function comprising up to $100$ segments. We select a one-dimensional domain $\Omega = [0, 2\pi]$, in which, given random variables $a_\ell$, $b_\ell$, and $c_\ell$, the subfunctions $a_\ell \exp(\sin(b_\ell x) + \cos(c_\ell x))$ are defined within each subdomain $\Omega_\ell = (\Gamma_\ell, \Gamma_{\ell+1})$. Here, the interfaces $\Gamma_\ell$ are uniformly distributed with $\Gamma_\ell = \ell\frac{2\pi}{L+1}$. Table~\ref{Table:1d_multi_piece_comparison} summarizes the approximation results for the cases with 5, 10, 50, and 100 subdomains. For the loss model, 1000 training points are sampled in the first three cases, while 2000 points are used for the 100-piece case.

\begin{table}[h!]
\caption{Numerical results for approximating a one-dimensional multi-piece function with different numbers of pieces in Example~2. CE: categorical embedding model; SE: scalar encoding model; OH: one-hot encoding model. ``--" means the method fails to converge.}
\label{Table:1d_multi_piece_comparison}
\centering
\begin{tabular}{c|c|ccc}
\hline
Number of pieces & Method & $N_p$ & $L^2$ error & $L^\infty$ error\\
\hline
\multirow{5}{*}{5} &CE ($D = 1$)                             & 205 & 5.65E$-$08 & 1.43E$-$07\\
                            &CE ($D = 2$)                             & 260 & 2.76E$-$08 & 8.14E$-$08\\
                            &SE ($\gamma_\ell = \ell$)  & 200 & 4.02E$-$08 & 3.43E$-$07\\
                            &SE ($\gamma_\ell = \bar{\gamma_\ell}$)  & 200 & 9.60E$-$08 & 2.19E$-$07\\
                            &OH                                             & 400 & 2.79E$-$08 & 7.14E$-$08\\
\hline
\multirow{6}{*}{10} &CE ($D = 1$)                              & 210 & 4.27E$-$08 & 3.96E$-$07\\
                              &CE ($D = 2$)                              & 270 & 4.78E$-$08 & 1.68E$-$07\\
                              &CE ($D = 5$)                              & 450 & 3.24E$-$08 & 1.87E$-$07\\
                              &SE ($\gamma_\ell = \ell$)  & 200 & -- & --\\
                              &SE ($\gamma_\ell = \bar{\gamma_\ell}$)  & 200 & 3.81E$-$08 & 5.67E$-$07\\
                              &OH                                             & 650 & 3.26E$-$08 & 1.49E$-$07\\
\hline
\multirow{6}{*}{50} &CE ($D = 1$)                              & 250 & 4.03E$-$04 & 2.94E$-$03\\
                              &CE ($D = 2$)                              & 350 & 1.71E$-$06 & 2.50E$-$05\\
                              &CE ($D = 5$)                              & 650 & 3.89E$-$07 & 5.29E$-$06\\
                              &CE ($D = 10$)                            & 1150 & 7.43E$-$08 & 1.47E$-$06\\
                              &SE ($\gamma_\ell = \bar{\gamma_\ell}$)  & 200 & 8.29E$-$04 & 8.92E$-$03\\
                              &OH                                             & 2650 & 4.29E$-$08 & 6.00E$-$07\\
\hline
\multirow{6}{*}{100} &CE ($D = 1$)                            & 300 & 2.27E$-$03 & 1.99E$-$02\\
                              &CE ($D = 2$)                              & 450 & 1.46E$-$05 & 1.92E$-$04\\
                              &CE ($D = 5$)                              & 900 & 1.94E$-$07 & 3.83E$-$06\\
                              &CE ($D = 10$)                            & 1650 & 9.02E$-$08 & 2.83E$-$06\\
                              &SE ($\gamma_\ell = \bar{\gamma_\ell}$)  & 200 & 4.47E$-$03 & 3.62E$-$02\\
                              &OH                                             & 5150 & 5.69E$-$08 & 1.43E$-$06\\
\hline
\end{tabular}
\end{table}

Notably, the one-hot model consistently achieves high prediction accuracy across all cases, succeeding even in the 100-piece case with $L^2$ error as low as $10^{-8}$. However, this model requires the largest number of trainable parameters, resulting in a substantial computational workload. By contrast, regardless of the number of segments, the design of the scalar encoding model with given labels ($\gamma_\ell$) requires the same (and least) number of trainable parameters $N_p = 200$. When using a simple nominal label $\gamma_\ell = \ell$, the scalar encoding model performs well only for the 5-piece case but fails in the other cases. To improve the model's capability, it is natural to deepen or widen the FC structure in the scalar encoding model. However, our experiments indicate that this strategy still fails to approximate discontinuous functions with many pieces (not shown here). This issue can be easily cured by setting the mean of the target label $\gamma_\ell = \bar{\gamma_\ell} = \int_{\Omega_\ell} u(\bx)\mbox{ d}\bx / |\Omega_\ell|$, which can be approximated simply via Monte Carlo integration. Despite the decrease in prediction accuracy with an increasing number of segments, this labeling strategy consistently enables successful training across all the cases, emphasizing the importance of incorporating informative categorical labels.

To investigate the effect of the reduced dimension $D$ in the categorical embedding model, we test various values of $D$ across all the cases. For the cases with 5 and 10 pieces, setting $D = 1$ is sufficient to achieve highly accurate prediction models with $L^2$ error as low as $10^{-8}$. Increasing $D$ to 2 or 5 provides only a minor improvement in accuracy. However, in the cases with 50 or 100 pieces, a higher-dimensional embedded space may be needed to capture more intrinsic features of the subfunctions. As a result, increasing $D$ generally leads to better approximation outcomes. Our experiments indicate that for cases with a large number of pieces, setting the reduced dimension $D$ to about 10\% or 20\% of the number of pieces yields a categorical embedding model with comparable accuracy to the one-hot encoding model, while requiring significantly fewer trainable parameters.

%%%%%%%%%%%%%%%%%%%%%%%%%%%%%%%%%%%%%%%%%
\subsection{Examples for anisotropic elliptic interface problems}\label{subsec:test_anisotropicEIP}
%%%%%%%%%%%%%%%%%%%%%%%%%%%%%%%%%%%%%%%%%
In this section, we present several examples for solving the anisotropic elliptic interface problems across one to three dimensions. To evaluate the performance of the proposed method, for each case, we derive the terms $f(\bx)$, $g(\bx)$, $v_\ell(\bx)$, and $w_\ell(\bx)$ from the exact solution $u(\bx)$, the coefficient matrix $A(\bx)$, and the scalar function $\lambda(\bx)$. With the knowledge of these terms, we then train the model to minimize the loss function (\ref{Eq:loss_aniso}). The training procedure is terminated when the loss value is smaller than a prescribed tolerance $\varepsilon = 10^{-15}$ or reaches 1000 training iterations. In Examples~1 to ~4, we deploy $N = 50$ neurons in the FC layer, while $N =100$ neurons for Example~5. Again, each neuron is activated by the sigmoid function. We also recall that the scalar encoding and one-hot encoding models are applied in the PINN-type loss~(\ref{Eq:loss_aniso}) simply by fixing $E = [\gamma_0, \gamma_1, \cdots, \gamma_{L}]^\top$ and $E = I_{L+1}$, respectively. Again, we report the average $L^2$ and $L^\infty$ errors using randomly sampled test points, which are 10 times the number of interior training points $M$, over 10 trial runs for each case.

%%%%%%%%%%%%%%%%%%%%%%%%%%%%%%%%%%%%%%%%%
\subsubsection*{\textbf{Example 1}}
%%%%%%%%%%%%%%%%%%%%%%%%%%%%%%%%%%%%%%%%%

In the first example, we demonstrate the capability of the proposed CE models for solving one-dimensional anisotropic problems with numerous jump discontinuities. Following the same setup as in Example~2 of Sec.~\ref{subsec:test_func_approx}, namely, we uniformly partition the domain $\Omega = [0, 2\pi]$ into subdomains $\Omega_\ell = (\Gamma_\ell, \Gamma_{\ell+1})$ with $\Gamma_\ell = \ell\frac{2\pi}{L+1}$. In each subdomain, we set the exact solution $u(x) = a_\ell \exp(\sin(b_\ell x) + \cos(c_\ell x))$, the anisotropic coefficient $A(x) = (d_\ell x)^2$, and $\lambda(x) = \sin^2(e_\ell x)$, with randomly predefined variables $(a_\ell, b_\ell, c_\ell, d_\ell, e_\ell)$. Notice that here we choose exactly the same random variables $a_\ell, b_\ell$ and $c_\ell$ as in the function approximation case by fixing the random seed.

Table~\ref{Table:1d_aniso} shows the results for cases with 5, 10, 50, and 100 pieces, exhibiting trends similar to those observed in the function approximation tests in Sec.~\ref {subsec:test_func_approx}. As expected, the one-hot encoding categorization model performs effectively across all cases, achieving accuracy with $L^2$ errors on the order of $10^{-8}$. However, this comes with a tradeoff between prediction accuracy and training cost due to the large number of trainable parameters. On the other hand, the scalar encoding with nominal labeling $\gamma_\ell = \ell$ succeeds only in the 5-piece case. As encountered in the function approximation experiments, this can be addressed by assigning a more informative mean value to $\gamma_\ell$. Here, we use the mean of the right-hand side function $f$ by setting $\bar{f_\ell} = \int_{\Omega_\ell} f(x)\, \text{d}x / |\Omega_\ell|$, and normalize those mean values to assign $\gamma_\ell = \bar{\gamma_\ell} = \bar{f_\ell}/\max_{1\leq\ell\leq L+1}|\bar{f_\ell}|$ in the embedding matrix since the $\bar{f_\ell}$ values range from $10^{-3}$ to $10^3$ in this experiment. This allows the scalar encoding model to solve the cases with 5 and 10 pieces, but remains insufficient for the cases with more than 50 pieces, leaving an open question regarding selecting appropriate labels $\gamma_\ell$.

\begin{table}[t]
\caption{Numerical results for solving the one-dimensional anisotropic elliptic interface problem with different numbers of pieces in Example~1. CE: categorical embedding model; SE: scalar encoding model; OH: one-hot encoding model. ``--" means the method fails to converge.}
\label{Table:1d_aniso}
\centering
\begin{tabular}{c|c|ccc}
\hline
Number of pieces & Method & $N_p$ & $L^2$ error & $L^\infty$ error\\
\hline
\multirow{5}{*}{5} &CE ($D = 1$)                             & 205 & 4.37E$-$07 & 1.78E$-$06\\
                            &CE ($D = 2$)                             & 260 & 6.13E$-$08 & 1.47E$-$07\\
                            &SE ($\gamma_\ell = \ell$)  & 200 & 7.63E$-$07 & 4.01E$-$06\\
                            &SE ($\gamma_\ell = \bar{\gamma_\ell}$)  & 200 & 6.40E$-$07 & 2.71E$-$06\\
                            &OH                                             & 400 & 2.57E$-$08 & 7.69E$-$08\\
\hline
\multirow{5}{*}{10} &CE ($D = 1$)                              & 210 & 1.68E$-$06 & 4.68E$-$06\\
                              &CE ($D = 2$)                              & 270 & 2.32E$-$07 & 6.43E$-$07\\
%                              &CE ($D = 5$)                              & 450 & 6.51E$-$08 & 1.80E$-$07\\
                              &SE ($\gamma_\ell = \ell$)  & 200 & -- & --\\
                              &SE ($\gamma_\ell = \bar{\gamma_\ell}$)  & 200 & 1.36E$-$03 & 3.25E$-$03\\
                              &OH                                             & 650 & 7.00E$-$08 & 2.20E$-$07\\
\hline
\multirow{4}{*}{50} &CE ($D = 5$)                              & 650 & 2.38E$-$05 & 1.44E$-$04\\
                              &CE ($D = 10$)                            & 1150 & 2.26E$-$05 & 9.29E$-$05\\
                              &SE ($\gamma_\ell = \bar{\gamma_\ell}$)  & 200 & -- & --\\
                              &OH                                             & 2650 & 3.47E$-$08 & 8.04E$-$07\\
\hline
\multirow{4}{*}{100} &CE ($D = 5$)                              & 900 & 5.77E$-$05 & 3.29E$-$04\\
                              &CE ($D = 10$)                            & 1650 & 4.97E$-$05 & 1.94E$-$04\\
                              &SE ($\gamma_\ell = \bar{\gamma_\ell}$)  & 200 & -- & --\\
                              &OH                                             & 5150 & 8.08E$-$08 & 4.23E$-$07\\
\hline
\end{tabular}
\end{table}

Setting the reduced dimension $D = 1$ or $2$ in the categorical embedding model performs effectively for the 5- and 10-piece cases, achieving accuracy comparable to the one-hot encoding model. However, unlike in the function approximation context, a low-dimensional embedded space ($D = 1$ or 2) in 50- and 100-piece cases may be inadequate to capture the intrinsic features of the network solution at the PDE level using the PINN learning machinary. As observed, increasing the reduced dimension $D$ up to 5 and 10 significantly enhances the model's capability, allowing it to learn the complex network solutions required for the 50- and 100-piece cases.

%%%%%%%%%%%%%%%%%%%%%%%%%%%%%%%%%%%%%%%%%
\subsubsection*{\textbf{Example 2}}
%%%%%%%%%%%%%%%%%%%%%%%%%%%%%%%%%%%%%%%%% 
Next, we turn to solve a two-dimensional problem with anisotropic variable coefficients, which often serves as a benchmark tested in various numerical methods, see~\cite{PWXY22, XZ24, DFL20}. The designated domain is a regular square $\Omega = [-1,1]^2$ containing a heart-shaped interface $\Gamma_1 = \{(r(\theta)\cos\theta-0.25,r(\theta)\sin\theta)|r(\theta) = (1+\cos\theta)/3, \theta\in[0,2\pi)\}$, which partitions the domain into two subdomains, $\Omega = \Omega_0\cup\Omega_1$. The exact solution is given by
\begin{equation*}
u(x_1,x_2) =
\begin{cases}
x_1^2+x_2^2            & \mbox{\;\;if\;\;} (x_1,x_2) \in \Omega_0,\\
\exp(x_1)\cos(x_2)  & \mbox{\;\;if\;\;} (x_1,x_2) \in \Omega_1.
\end{cases}
\end{equation*}
Define
\begin{eqnarray*}
A_1(x_1,x_2) =
\left[
\begin{array}{cc}
x_1^2 + x_2^2+1 & x_1^2 + x_2^2\\
x_1^2 + x_2^2 & x_1^2 + x_2^2+2
\end{array}\right],
\end{eqnarray*}
and $\lambda_1(x_1,x_2) = \exp(x_1)(x_1^2 + x_2^2+3)\sin(x_2)$,
the anisotropic coefficient $A$ and the scalar function $\lambda$ are piecewise spatial dependent functions set by
\begin{eqnarray*}
A(x_1,x_2) =
\begin{cases}
1000A_1(x_1,x_2)& \mbox{\;\;if\;\;} (x_1,x_2) \in \Omega_0,\\
A_1(x_1,x_2)     & \mbox{\;\;if\;\;} (x_1,x_2) \in \Omega_1,
\end{cases}
\end{eqnarray*}
\begin{eqnarray*}
\lambda(x_1,x_2) =
\begin{cases}
1000\lambda_1(x_1,x_2)    & \mbox{\;\;if\;\;} (x_1,x_2) \in \Omega_0,\\
\lambda_1(x_1,x_2) & \mbox{\;\;if\;\;} (x_1,x_2) \in \Omega_1,
\end{cases}
\end{eqnarray*}
so that the contrasts for both functions are $1000$. 

In this test, we train the loss model using randomly selected training points with $(M, M_b, M_{\Gamma_1}) = (324, 72, 72)$. The results, compared with those from the finite volume method (FVM)~\cite{PWXY22}, are presented in Table~\ref{Table:ansio_cusp}. Notice that for FVM, the total number of degrees of freedom (or unknowns) is based on the number of discretization grid points, $m^2$. As shown in Table~\ref{Table:ansio_cusp}, FVM uses a grid resolution of \(m = 128\), resulting in \(N_{p} = 16384\) unknowns, while the neural network models require only a few hundred parameters. As noted, with a relatively small training dataset (a few hundred points) and a moderate number of trainable parameters (also in the few hundred range), all neural network models achieve higher accuracy than FVM. Given that the solution involves only two subdomain solutions, all models are capable of achieving similar levels of prediction accuracy.

\begin{table}[h!]
\caption{Numerical results for solving the two-dimensional anisotropic interface problem in Example~2. CE: categorical embedding model; SE: scalar encoding model; OH: one-hot encoding model; FVM: finite volume method~\cite{PWXY22}.}
\label{Table:ansio_cusp}
\centering
\begin{tabular}{c|ccc}
\hline
Method & $N_p$ & $L^2$ error & $L^\infty$ error\\
\hline
CE ($D = 1$) & 252 & 5.97E$-$08 & 5.43E$-$07\\
SE ($\gamma_\ell = \ell$)         & 250 & 6.62E$-$08 & 6.23E$-$07\\
OH                & 300 & 1.47E$-$08 & 1.34E$-$07\\
FVM~\cite{PWXY22} & 16384 & 1.80E$-$04 \\
\hline
\end{tabular}
\end{table}

\begin{figure*}[h]
\centering
\includegraphics[scale=0.45]{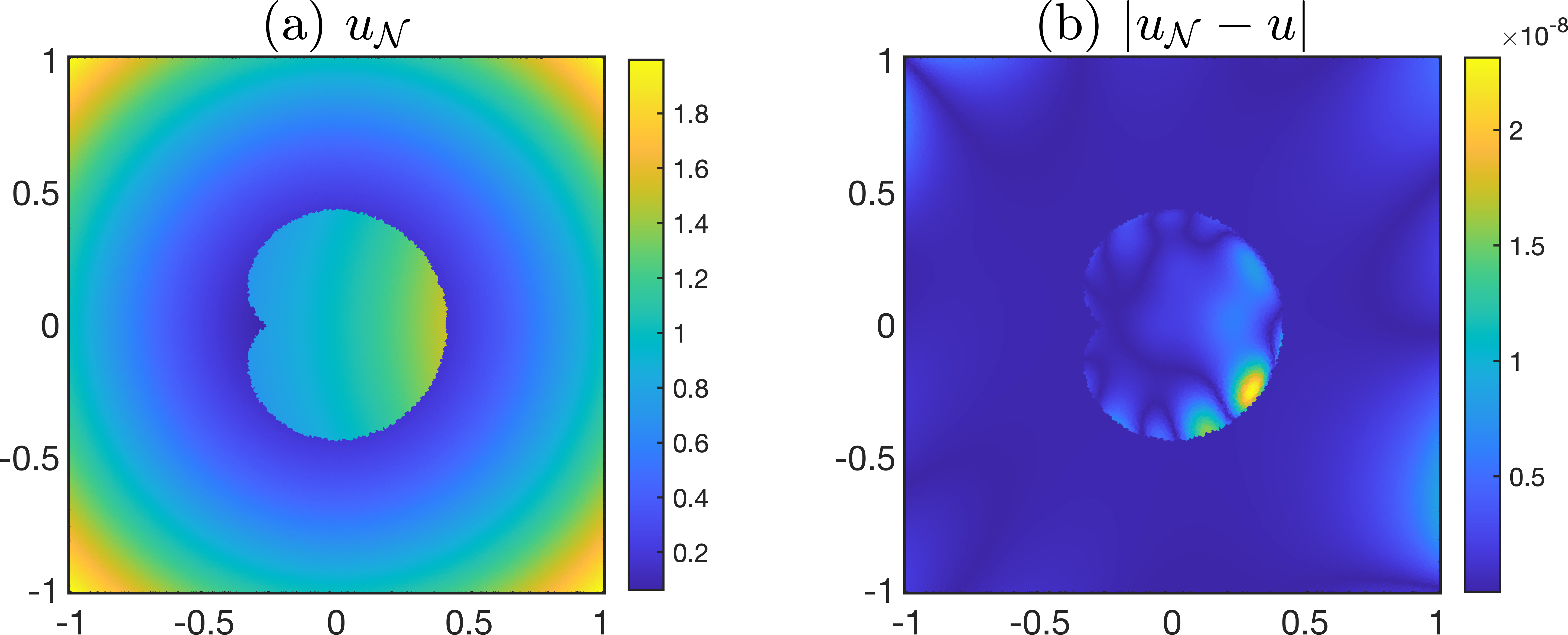}
\caption{(a) The trained CE model $u_\mathcal{N}$ for the piecewise-defined solution of Example~2. (b) The pointwise absolute error $|u_\mathcal{N}-u|$. The maximum error is $\|u_\mathcal{N}-u\|_\infty = 2.31\times10^{-8}$.}
\label{Fig:aniso_cusp}
\end{figure*}

The categorical embedding solution with reduced dimension \(D = 1\) is depicted in Fig.~\ref{Fig:aniso_cusp}(a). As shown, the model captures the jump discontinuity sharply along the interface, and the absolute error, depicted in Fig.~\ref{Fig:aniso_cusp}(b), exhibits a pointwise absolute error as low as \(10^{-8}\), demonstrating the high prediction accuracy of our proposed model. It is worth noting that, in this example, a cusp occurs at $(-0.25, 0)$, our discontinuity capturing models can handle such interfaces without any difficulty. In contrast, some existing numerical methods may encounter challenges with these singular points.

%%%%%%%%%%%%%%%%%%%%%%%%%%%%%%%%%%%%%%%%%
\subsubsection*{\textbf{Example 3}}
\begin{figure*}[h]
\centering
\includegraphics[scale=0.45]{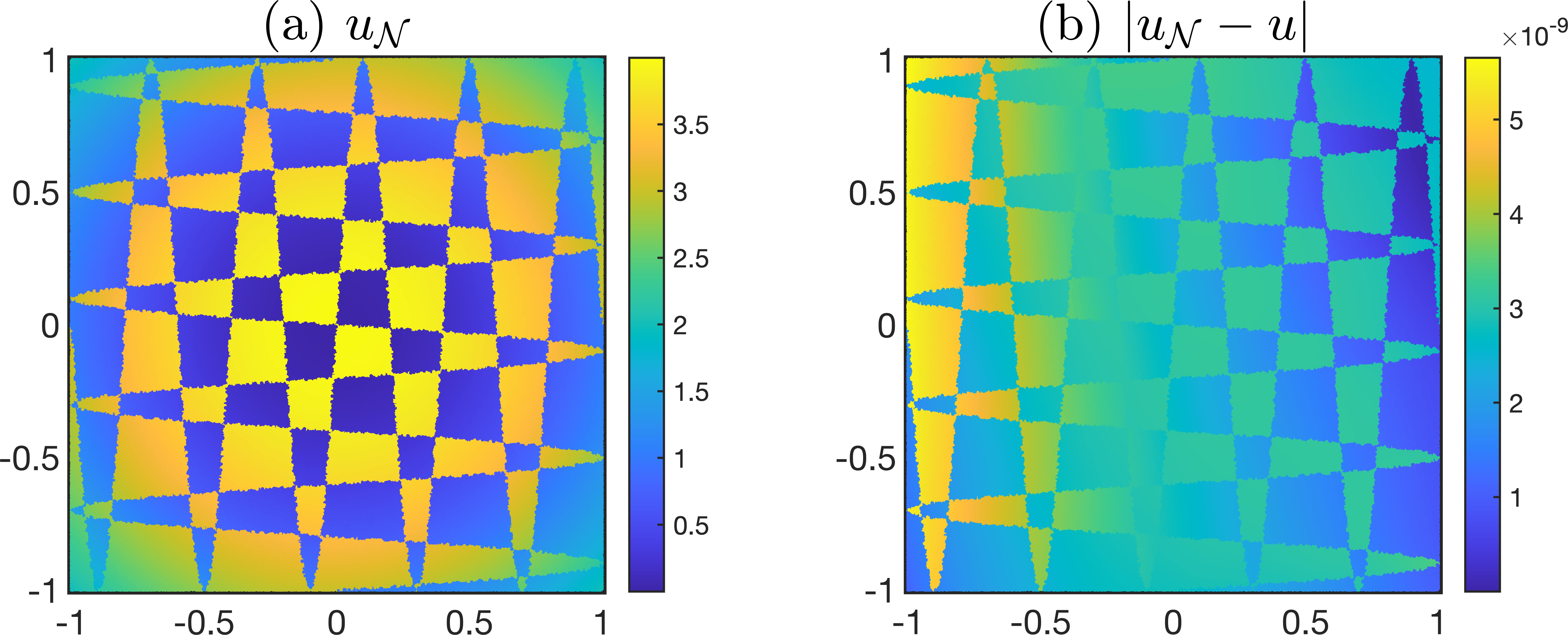}
\caption{(a) The trained categorical embedding solution $u_\mathcal{N}$ with $D=1$ of Example 3. (b) The pointwise absolute error $|u_\mathcal{N}-u|$. The maximum error is $\|u_\mathcal{N}-u\|_\infty = 5.65\times10^{-9}$.}
\label{Fig:aniso_chessboard}
\end{figure*}

%%%%%%%%%%%%%%%%%%%%%%%%%%%%%%%%%%%%%%%%%
In this example, we consider a benchmark for a 2D variable coefficient elliptic interface problem, see Example 4 in~\cite{HWW10}. With a square domain, $\Omega=[-1,1]^2$, the interface is described as the longitude and latitude of a ``chessboard''-like domain (see the left panel of Fig.~\ref{Fig:aniso_chessboard}) given by the zero level set of $\phi(x_1,x_2) = (\sin(5\pi x_1)-x_2)(-\sin(5\pi x_2)-x_1)$. Thus we define the separated domain $\Omega_0 = \{(x_1,x_2)|\phi(x_1,x_2)>0\}$ and $\Omega_1 = \{(x_1,x_2)|\phi(x_1,x_2)<0\}$. The diffusion coefficient is given by
\begin{equation*}
A(x_1, x_2) =
\begin{cases}
x_1x_2+2    & \mbox{\;\;if\;\;} (x_1, x_2) \in \Omega_0,\\
x_1^2-x_2^2+3 & \mbox{\;\;if\;\;} (x_1, x_2) \in \Omega_1.
\end{cases}
\end{equation*}
We set $\lambda=0$ and the exact solution
\begin{equation*}
u(x_1, x_2) =
\begin{cases}
4-x_1^2-x_2^2    & \mbox{\;\;if\;\;} (x_1, x_2) \in \Omega_0,\\
x_1^2+x_2^2 & \mbox{\;\;if\;\;} (x_1, x_2) \in \Omega_1.
\end{cases}
\end{equation*}
We follow the same setup as in Example~2. The categorical embedding prediction solution with $D=1$ along with its pointwise absolute error is displayed in Fig.~\ref{Fig:aniso_chessboard}. The present model can achieve a very accurate result with the $L^\infty$ error, $\|u_\mathcal{N}-u\|_\infty = 5.65\times10^{-9}$.

We summarize the accuracy between the present network models and the one of finite element method~(FEM)~\cite{HWW10} in Table~\ref{Table:ansio_chessboard}. Indeed, the prediction accuracy of all models clearly outperforms the one obtained by FEM; each network model with a few hundred trainable parameters readily achieves high solution expressivity with accuracy of order $10^{-9}$ in $L^2$ error.
It is worth noting that, although this problem can be dealt appropriately with traditional numerical methods such as FVM or FEM, these methods require identifying regular and irregular grid points (or cell triangulations) as a preliminary step, making their implementation somewhat tedious. On the contrary, it is straightforward to implement the present network models simply using the categorization map $\bz$. Moreover, we randomly sampled $m_{\Gamma_1} = 72$ points to represent the interface, and the optimization algorithm's implementation requires no additional effort compared to Examples~1 and 2, demonstrating the robustness of the present method regardless of the complexity of embedded interface geometries.

\begin{table}[h!]
\caption{Numerical results for solving the two-dimensional anisotropic interface problem in Example~3. CE: categorical embedding model; SE: scalar encoding model; OH: one-hot encoding model; FEM: finite element method~\cite{HWW10}.}
\label{Table:ansio_chessboard}
\centering
\begin{tabular}{c|ccc}
\hline
Method & $N_p$ & $L^2$ error & $L^\infty$ error\\
\hline
CE ($D = 1$) & 252 & 4.39E$-$09 & 9.59E$-$09\\
SE ($\gamma_\ell = \ell$)         & 250 & 5.48E$-$09 & 1.25E$-$08\\
OH                & 300 & 4.13E$-$09 & 7.89E$-$09\\
FEM~\cite{HWW10} & 102400 & & 2.60E$-$05  \\
\hline
\end{tabular}
\end{table}

%%%%%%%%%%%%%%%%%%%%%%%%%%%%%%%%%%%%%%%%%
\subsubsection*{\textbf{Example 4}}
%%%%%%%%%%%%%%%%%%%%%%%%%%%%%%%%%%%%%%%%%

In this example, we aim to highlight the capability of the present method by tackling two-dimensional problems with multiple subdomains enclosed in an irregular domain. Here, the domain and subdomain geometries are shown in Fig.~\ref{Fig:domain_2d}, while the detailed formulation can be found in Example~1 of Sec.~\ref{subsec:test_func_approx}. The solution profile is given in Eq.~(\ref{Eq:function_approx}); the coefficient matrix $A$ and scalar function $\lambda$ in $\Omega_0$ are respectively set by
\[
A_0(x_1,x_2) =
\left[
\begin{array}{cc}
(x_1+x_2)^2+1 & -x_1^2+x_2^2\\
-x_1^2+x_2^2 & (x_1-x_2)^2+1
\end{array}\right] \quad \mbox{and} \quad
\lambda_0(x_1,x_2) = \exp(x_1-x_2),
\]
while in the other subdomains, we set $A_\ell(x_1,x_2) = \beta_\ell A_0(x_1,x_2)$ and $\lambda_\ell(x_1,x_2) = \beta_\ell\lambda_0(x_1,x_2)$, where $\beta_1 = 10^{-1}$, $\beta_2 = 10^{-2}$, $\beta_3 = 10^{1}$, and $\beta_4 = 10^{2}$ (so the largest ratio in this case is $10000$). 

Following the same training setup as in Example~2, each solution model achieves high prediction accuracy, with \(L^2\) errors as low as \(10^{-9}\), as shown in Table~\ref{Table:ansio_2d_multipiece}. Notably, at the PDE-solving level, the categorical embedding model requires only about half the number of trainable parameters compared to the one-hot encoding model. Additionally, owing to the mesh-free nature of the neural network method, implementing the model is straightforward, as demonstrated in this test with the superellipse. In contrast, grid-based methods require substantially greater effort in implementation to solve problems on irregular domains with multiple subdomains.

\begin{table}[h!]
\caption{Numerical results for solving the two-dimensional anisotropic interface problem in Example 4. CE: categorical embedding model; SE: scalar encoding model; OH: one-hot encoding model.}
\label{Table:ansio_2d_multipiece}
\centering
\begin{tabular}{c|ccc}
\hline
Method & $N_p$ & $L^2$ error & $L^\infty$ error\\
\hline
CE ($D = 1$) & 255 & 2.33E$-$09 & 2.28E$-$08\\
SE ($\gamma_\ell = \ell$)         & 250 & 3.96E$-$09 & 3.54E$-$08\\
OH                & 450 & 2.95E$-$09 & 2.09E$-$08\\
\hline
\end{tabular}
\end{table}

%%%%%%%%%%%%%%%%%%%%%%%%%%%%%%%%%%%%%%%%%
\subsubsection*{\textbf{Example 5}}
%%%%%%%%%%%%%%%%%%%%%%%%%%%%%%%%%%%%%%%%%

\begin{figure}[h]
\centering
\includegraphics[scale=0.4]{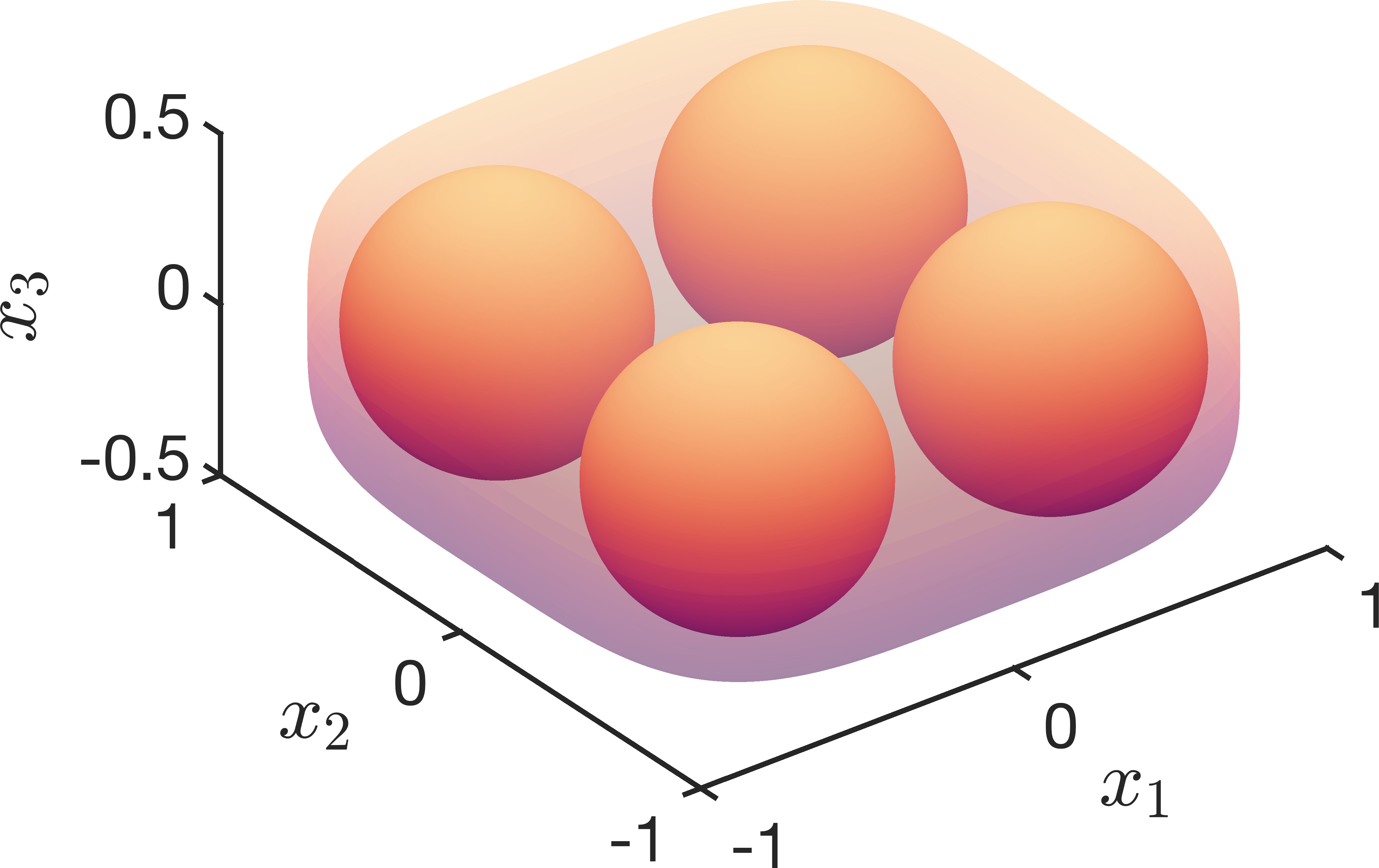}
\caption{The super-quadric domain $x_1^4+x_2^4+16x_3^4 = 1$ with four embedded spheres of radius $0.4$ located at $(-0.45,0.45,0)$, $(0.45,0.45,0)$, $(-0.45,-0.45,0)$, and $(0.45,-0.45,0)$.}
\label{Fig:surface_mp}
\end{figure}

In the last example, we illustrate the robustness of our method for solving a three-dimensional anisotropic problem with multiple interfaces. We set the domain $\Omega$ with a super-quadric boundary $x_1^4+x_2^4+16x_3^4 = 1$, in which there are four subdomains, $\Omega_1$, $\Omega_2$, $\Omega_3$, and $\Omega_4$, encapsulated by four spheres of radius $0.4$ with their center respectively located at $(-0.45,0.45,0)$, $(0.45,0.45,0)$, $(-0.45,-0.45,0)$, and $(0.45,-0.45,0)$. See the domain depiction in Fig.~\ref{Fig:surface_mp}.

Denoting $\bx = (x_1, x_2, x_3)$, we set the solution as
\begin{equation}
u(\bx) =
\begin{cases}
\exp(x_1+x_2+x_3) & \mbox{\;\;if\;\;} \bx \in \Omega_0,\\
\sin x_1\sin x_2\sin x_3 & \mbox{\;\;if\;\;} \bx \in \Omega_1,\\
\cos x_1\cos x_2\cos x_3 & \mbox{\;\;if\;\;} \bx \in \Omega_2,\\
\sinh x_1\sinh x_2\sinh x_3 & \mbox{\;\;if\;\;} \bx \in \Omega_3,\\
\cosh x_1\cosh x_2\cosh x_3 & \mbox{\;\;if\;\;} \bx \in \Omega_4,
\end{cases}
\end{equation}
and choose $A_0(\bx) = R\Lambda R^T$ and $\lambda_0(\bx) = \exp(x_1-x_2-x_3)$, where
\begin{equation*}
R =
\left[
\begin{array}{rrr}
2/3 & 1/3 & 2/3\\
-2/3 & 2/3 & 1/3\\
1/3 & 2/3 & -2/3
\end{array}\right] \quad \mbox{and} \quad
\Lambda =
\left[
\begin{array}{ccc}
\|\bx\|^2+1 & 0 & 0\\
0 & \|\bx\|^2+2 & 0\\
0 & 0 & \|\bx\|^2+3
\end{array}\right].
\end{equation*}
We further define \(A_\ell = \beta_\ell A_0\) and \(\lambda_\ell = \beta_\ell \lambda_0\), where \(\beta_1 = 0.1\), \(\beta_2 = 0.05\), \(\beta_3 = 10\), and \(\beta_4 = 50\), resulting in a maximum contrast ratio of \(1000\). In Table~\ref{Table:ansio_3d_multipiece}, we train each model using $(M,M_b,M_\Gamma) = (324, 144, 144)$.
In this setup, the categorical embedding network with \(D=1\) and the one-hot encoding model exhibit similar performance, achieving approximately one order of magnitude higher prediction accuracy than the scalar encoding model. Notably, the categorical embedding model learns the categorization map using a single-dimensional representation. This results in roughly half the number of trainable parameters compared to the one-hot encoding model, significantly reducing computational effort.

\begin{table}[h!]
\caption{Numerical results for solving the three-dimensional anisotropic interface problem in Example 5. CE: categorical embedding model; SE: scalar encoding model; OH: one-hot encoding model.}
\label{Table:ansio_3d_multipiece}
\centering
\begin{tabular}{c|ccc}
\hline
Method & $N_p$ & $L^2$ error & $L^\infty$ error\\
\hline
CE ($D = 1$) & 605 & 3.42E$-$08 & 2.64E$-$07\\
SE ($\gamma_\ell = \ell$)         & 600 & 2.87E$-$07 & 3.45E$-$06\\
OH                & 1000 & 4.25E$-$08 & 5.51E$-$07\\
\hline
\end{tabular}
\end{table}

%%%%%%%%%%%%%%%%%%%%%%%%%%%%%%%%%%%%%%%%%
\section{Conclusion}
%%%%%%%%%%%%%%%%%%%%%%%%%%%%%%%%%%%%%%%%%
In this work, we propose a single neural network architecture, a categorical embedding discontinuity-capturing shallow neural network, for representing piecewise smooth functions. The architecture consists of three hidden layers: (i) a discontinuity-capturing layer, which maps domain segments to disconnected sets in a higher-dimensional space; (ii) a categorical embedding layer, which reduces the high-dimensional information into low-dimensional features; (iii) a fully connected layer, which models the continuous mapping. This design enables the accurate approximation of piecewise smooth functions, even when a function contains a large number of pieces. Moreover, the proposed network possesses two notable features: it allows for the direct calculation of jump values at interfaces, and it ensures well-defined derivatives at all points away from interfaces and domain boundaries, thereby facilitating derivative evaluations.

We further apply the proposed categorical embedding discontinuity capturing shallow neural network to solve anisotropic elliptic interface problems, which are traditionally regarded as highly challenging. Follow the PINN framework, the network is trained using the LM optimizer by minimizing the mean squared error loss of the system. With such a simple design and shallow network architecture, the model achieves accuracy and efficiency comparable to established grid-based numerical methods. Importantly, the approach is entirely mesh-free: once training points are appropriately selected, the loss formulation and optimization process remain unchanged across different domains and subdomain geometries.

Finally, our results demonstrate that machine learning-based methods can achieve accuracy comparable to conventional scientific computing approaches, while offering additional advantages. These include their mesh-free nature, ease of implementation, natural compatibility with GPU acceleration, and flexibility for extension to high-dimensional problems. Such properties suggest strong potential for practical applications in computational physics and engineering, where both accuracy and efficiency are critical.

%%%%%%%%%%%%%%%%%%%%%%%%%%%%%%%%%%%%%%%%%
\section*{Acknowledgments}
%%%%%%%%%%%%%%%%%%%%%%%%%%%%%%%%%%%%%%%%%

W.-F. Hu, T.-S. Lin, Y.-H. Tseng, and M.-C. Lai acknowledge the supports by National Science and Technology Council, Taiwan, under research grants 114-2628-M-008-002-MY4, 111-2628-M-A49-008-MY4, 113-2115-M-390-007-MY2, and 113-2115-M-A49-014-MY3, respectively. W.-F. Hu and T.-S. Lin also acknowledge the supports by National Center for Theoretical Sciences, Taiwan.

\end{document}